\documentclass[preprint,12pt]{elsarticle}
\usepackage[utf8]{inputenc}
\usepackage[margin=2.54cm]{geometry}

\usepackage{amsmath}
\usepackage{amssymb}
\usepackage{amsthm}
\usepackage{bbm}
\usepackage{xcolor}
\usepackage{hyperref}
\usepackage{mathtools}
\usepackage{enumitem}
\usepackage{graphicx}
\usepackage{caption}
\usepackage{subcaption}
\usepackage{changepage}

\DeclareMathOperator*{\expct}{\mathbb{E}}

\makeatletter
    \newtheorem*{rep@theorem}{\rep@title}
    \newcommand{\newreptheorem}[2]{
    \newenvironment{rep#1}[1]{
     \def\rep@title{#2 \ref{##1}}
     \begin{rep@theorem}}
     {\end{rep@theorem}}}
\makeatother

\newcommand{\vecsp}[2]{\mathbb{F}_{#1}^{#2}}
\newcommand{\id}[1]{\mathbbm{1}_{#1}}
\newcommand{\st}{:\mathrm{ }}
\newcommand{\twr}{\mathrm{\textbf{twr}}}
\newcommand{\wwz}{\mathrm{\textbf{wwz}}}
\newcommand{\bracketed}[1]{{(#1)}}
\newcommand{\vecline}[1]{\underline{#1}}
\newcommand{\rank}{\textbf{rank}}

\newcommand{\cmnt}[1]{\textcolor{black}{#1}}
\newcommand{\Img}[1]{\mathrm{Im}(#1)}

\newcommand{\codim}{\mathrm{codim}}

\DeclarePairedDelimiter\gen{\langle}{\rangle}
\DeclarePairedDelimiter\abs{|}{|}
\DeclarePairedDelimiter\norm{\|}{\|}
\DeclarePairedDelimiter\babs{\big|}{\big|}
\DeclarePairedDelimiter\set{\{}{\}}

\DeclarePairedDelimiter\ceil{\lceil}{\rceil}

\begin{document}

\newtheorem{theorem}{Theorem}[section]
\newtheorem{definition}[theorem]{Definition}
\newtheorem{conjecture}[theorem]{Conjecture}
\newtheorem*{arl}{Arithmetic Regularity Lemma}
\newtheorem*{nntheorem}{Theorem}
\newreptheorem{lemma}{Lemma}
\newtheorem{lemma}[theorem]{Lemma}
\newtheorem*{nclaim}{Claim}
\newtheorem{claim}[theorem]{Claim}
\newtheorem{corollary}[theorem]{Corollary}
\newtheorem{prop}[theorem]{Proposition}
\newtheorem{conj}[theorem]{Conjecture}

\theoremstyle{definition}
\newtheorem{example}[theorem]{Example}

\begin{frontmatter}

    \title{Induced arithmetic removal for partition-regular patterns of complexity 1}
    
    \author{V. Gladkova}
    \affiliation{organization={DPMMS, University of Cambridge},
                addressline={Centre for Mathematical Sciences}, 
                city={Cambridge},
                postcode={CB3 0WB},
                country={UK}}
    
    \begin{abstract}
        In 2019, Fox, Tidor and Zhao \cite{induced-1} proved an induced arithmetic removal lemma  for linear patterns of complexity 1 in vector spaces over a fixed finite field. With no further assumptions on the pattern, this induced removal lemma cannot guarantee a fully pattern-free recolouring of the space, as some `non-generic' instances must necessarily remain. On the other hand, Bhattacharyya, Fischer, H.~Hatami, P.~Hatami, and Lovett \cite{trans-invariant} showed in 2012 that in the case of translation-invariant patterns, it is possible to obtain recolourings that eliminate the given pattern completely, with no exceptions left behind. This paper demonstrates that such complete removal can be achieved for all partition-regular patterns of complexity 1.
    \end{abstract}
    
    %% Keywords

\end{frontmatter}

\section{Introduction}
The triangle removal lemma for graphs and, more generally, arbitrary subgraph removal lemmas \cite{removal-lemmas} are a well-known application of Szemer\'edi's regularity lemma \cite{szemeredi}. These removal results state that if a graph $G$ contains few copies of a graph $H$ as a subgraph, then $G$ can be made $H$-free by flipping only a small proportion of its edges. The same holds for $H$ occuring as an induced subgraph, as proved by Alon, Fischer, Krivelevich, and Szegedy~\cite{induced-graph} via a strong variant of the regularity lemma, which produces two nested vertex partitions with certain desirable properties.

Analogous induced removal lemmas in the setting of vector spaces over finite fields assert that if $\vecsp{p}{n}$ is $r$-coloured and a given arithmetic pattern occurs with small density under this colouring, then it is possible to recolour a small proportion of the space to eliminate (nearly) all instances of the pattern. Here, an arithmetic pattern is a coloured solution to a given system of equations with a specified assignment of colours. The following definition specifies a pattern using a collection of linear forms but this can be readily translated to the language of systems of equations by standard methods, as laid out in Section \ref{section:patterns}.

\begin{definition}[Arithmetic patterns]
\label{def:pattern}
An \emph{arithmetic pattern} $\mathcal{H}$ in $\ell$ variables is a tuple $(\mathcal{L}, \mathcal{X})$, where $\mathcal{L}$ is a collection of linear forms $\set{L_1, \ldots, L_m}$ in $\ell$ variables over $\vecsp{p}{}$, and $\mathcal{X}=\set{\chi:[m] \rightarrow [r]}$ is a collection of $r$-colourings.

Given an $r$-colouring $\phi$ of $\vecsp{p}{n}$, an \emph{instance of $\mathcal{H}$ under $\phi$} is an $\ell$-tuple $\vecline{x} \in (\vecsp{p}{n})^{\cmnt{\ell}}$ such that for some $\chi \in \mathcal{X}$, $\phi(L_j(\vecline{x})) = \chi(j) \text{ for each }1 \leq j \leq m.$ A colouring $\phi$ is said to be \emph{$\mathcal{H}$-free} if there are no instances of $\mathcal{H}$ under $\phi$.
\end{definition}

\noindent \textbf{Note.} Unlike in previous work, this definition of an arithmetic pattern incorporates multiple colourings at once. In this way, a monochromatic 3-term arithmetic progression can be expressed as a single pattern instead of a union of $r$ separate patterns. For example, when $r = 2$, this is given by $(\set{x, x+d, x+2d}, \set{\chi_1, \chi_2})$, where $\chi_i: [3] \rightarrow [2]$ is the constant colouring $\chi_i \equiv i$.

Induced arithmetic removal lemmas were established in the work of \cite{induced-over-f2, trans-invariant, induced-1}, and \cite{full-induced}, but, unlike in the graph-theoretic setting, the results of both \cite{induced-1} and \cite{full-induced} only hold with the caveat that a small number of `non-generic' instances might remain behind. In the case of complexity 1 patterns, Fox, Tidor and Zhao \cite{induced-1} proved the following, with $\Lambda_{\mathcal{H}}(\phi)$ denoting the density of the instances of $\mathcal{H}$ under a colouring $\phi$ (see Definition \ref{def:pattern-density}).

\begin{theorem}[Induced removal for complexity 1 patterns \cite{induced-1}]
    \label{thm:induced-removal-general-1}
    Fix $\epsilon > 0$, an integer $r > 0$, and an arithmetic pattern $\mathcal{H}$ of complexity 1. There exists a $\delta = \delta(\epsilon, r, \mathcal{H})$ satisfying the following. If $\phi:\vecsp{p}{n} \rightarrow [r]$ is an $r$-colouring of $\vecsp{p}{n}$ such that $\Lambda_{\mathcal{H}}(\phi) \leq \delta$, then $\phi$ can be made $\mathcal{H}$-free on $\vecsp{p}{n}\backslash\set{0}$ by recolouring at most an $\epsilon$-proportion of $\vecsp{p}{n}$.
\end{theorem}

\noindent This result cannot be extended to recolour all of $\vecsp{p}{n}$ when the arithmetic pattern in question is not partition-regular (see Example \ref{non-example}). The definition of partition regularity in the context of vector spaces over finite fields follows the work of Bergelson, Deuber, and Hindman \cite{partion-regular-definition}.
\begin{definition}[Partition regularity]
    \label{def:partition-regular}
    A linear system $\mathcal{L} = \set{L_1, \ldots, L_m}$ over $\vecsp{p}{}$ is said to be \emph{partition-regular} if the following holds. For any $r \in \mathbb{N}$, there exists $N = N(r, p, m)$ such that, if $n \geq N$, there is a monochromatic instance of $\mathcal{L}$ in $\vecsp{p}{n}\backslash\set{0}$ under any $r$-colouring of the space. An arithmetic pattern $\mathcal{H} = (\mathcal{L}, \mathcal{X})$ is partition-regular if the underlying linear system $\mathcal{L}$ is partition-regular.
\end{definition}

\begin{example}
\label{non-example}
Consider any linear system $\mathcal{L}$ over $\vecsp{p}{}$ that is not partition-regular. Then for any $r, n \in \mathbb{N}$, there is an $r$-colouring of $\vecsp{p}{n}$ such that the only monochromatic \cmnt{instance of $\mathcal{L}$} is $0$. However, this remains a monochromatic \cmnt{instance} no matter how the space is recoloured. (A special case of this was given in \cite[Non-example 1.3]{induced-1}.)
\end{example}

\noindent As a consequence, any induced arithmetic removal lemma for general arithmetic patterns must necessarily make exceptions for instances of a certain form. However, for partition-regular patterns (including translation-invariant patterns), an obstacle such as in Example \ref{non-example} does not apply: indeed, any colouring of $\vecsp{p}{n}$ must have \textit{many} monochromatic instances of such patterns (see Theorem \ref{theorem:rado}), whereas the induced removal lemma only applies when the number of instances is small to begin with.

In fact, the induced removal lemma of Bhattacharyya et al.~\cite{trans-invariant} for translation-invariant patterns guarantees the removal of \textit{all} instances. Such complete removal is possible due to the additional flexibility afforded by being able to translate instances of the pattern without breaking their structure (see the discussion in Section \ref{section:argument-overview}).

The main result of this paper demonstrates that it is possible to arrange for similar flexibility and therefore achieve complete removal in the case of partition-regular patterns of complexity 1.

\begin{theorem}[Induced removal for partition-regular patterns]
\label{main}
Fix $\epsilon > 0$, an integer $r > 0$, and a partition-regular pattern $\mathcal{H}$ of complexity 1. There exists a $\delta = \delta(\epsilon, r, \mathcal{H})$ satisfying the following. If the set of instances of $\mathcal{H}$ has density at most $\delta$ under an $r$-colouring $\phi:\vecsp{p}{n} \rightarrow [r]$, then $\phi$ can be made $\mathcal{H}$-free by recolouring at most an $\epsilon$-proportion of $\vecsp{p}{n}$.
\end{theorem}

\noindent The following example serves as an illustration of this theorem.

\begin{example}
    Let $\mathcal{H}$ be the arithmetic pattern encoding rainbow solutions to $x+y+z-w=0$ under $\phi$, noting that this is a partition-regular pattern of complexity 1 which is not translation-invariant when $p > 2$. For a given $\epsilon > 0$, let $\phi: \vecsp{p}{n} \rightarrow [4]$ be any surjective 4-colouring such that $0 \in \phi^{-1}(1)$ and $\abs{\phi^{-1}(1)} \leq \delta \abs{\vecsp{p}{n}}/4$, where $\delta = \delta_{\ref{main}}(\epsilon, 4, \mathcal{H})$. Since one of the colour classes is small, there are at most $\delta\abs{\vecsp{p}{n}}^3$ instances of $\mathcal{H}$ under $\phi$. Then Theorem \ref{main} implies that it is possible to eliminate all rainbow solutions to $x+y+z-w=0$ by recolouring at most an $\epsilon$-proportion of $\vecsp{p}{n}$.
    
    Of course, in this example, this is easy to see without resorting to Theorem \ref{main}, as we can simply replace all of $\phi^{-1}(1)$ with another colour. However, note that the induced removal lemma of Fox, Tidor and Zhao \cite{induced-1} (Theorem \ref{thm:induced-removal-general-1}) applied to this same setting would only eliminate rainbow solutions for which each of $x, y, z, w$ is non-zero.
\end{example}

\noindent Additionally, an easy corollary is provided here as another instance of the kind of results implied by Theorem \ref{main}. Here $A \Delta H$ denotes the symmetric difference of $A$ and $H$, so that \ref{epsilon-close-to-H} says that $A$ is  $\epsilon$-close to being a subspace.

\begin{corollary}
    Fix $\epsilon > 0$, and let $A$ be a subset of $\vecsp{p}{n}$. There exists a $\delta = \delta(\epsilon)$ such that one of the following holds:
    \begin{enumerate}[label = (\roman*)]
        \item \label{A-is-small} $\abs{A} \leq \epsilon \abs{\vecsp{p}{n}}$;
        \item $\abs{(A+A) \cap A^c} > \delta \abs{\vecsp{p}{n}}$;
        \item \label{epsilon-close-to-H} there is a subspace $H \leqslant \vecsp{p}{n}$ such that $\abs{A \Delta H} \leq \epsilon \abs{\vecsp{p}{n}}$.
    \end{enumerate}
\end{corollary}
\begin{proof}
    Let $\phi:\vecsp{p}{n} \rightarrow \set{0,1}$ be the 2-colouring $\phi = \id{A}$, and let $\mathcal{H}$ be the arithmetic pattern corresponding to the solutions of $x+y=z$ for which $x, y$ have colour $1$ and $z$ has colour $0$. This is a partition-regular pattern of complexity 1, so Theorem \ref{main} applies. If $\abs{(A+A) \cap A^c} \leq \delta \abs{\vecsp{p}{n}}$ where $\delta = \delta_{\ref{main}}(\epsilon, 2, \mathcal{H})$, then there are at most $\delta\abs{\vecsp{p}{n}}^2$ such solutions to $x+y=z$ under $\phi$. Therefore, by Theorem \ref{main}, there is a 2-colouring $\psi:\vecsp{p}{n} \rightarrow \set{0,1}$ that differs from $\phi$ in at most $\epsilon \abs{\vecsp{p}{n}}$ places such that $\psi$ has no solutions to $x+y=z$ satisfying $\psi(x)=\psi(y)=1$ and $\psi(z)=0$.

    In particular, letting $A' = \psi^{-1}(1)$, this implies that $A'+A' \subseteq A'$. Hence $A'$ is either a subspace or an empty set. Moreover, $\abs{A \Delta A'} \leq \epsilon \abs{\vecsp{p}{n}}$, as $A \Delta A'$ is precisely equal to the set of $x \in \vecsp{p}{n}$ such that $\phi(x) \neq \psi(x)$. Therefore either \ref{A-is-small} or \ref{epsilon-close-to-H} holds, as required.
\end{proof}

The proof of Theorem \ref{main} follows the same lines as existing induced removal lemmas but utilises a more general subcoset selection scheme, reminiscent of the argument for the induced graph removal lemma \cite{induced-graph}. The proposed approach is sketched out in Section \ref{section:argument-overview} before being formally applied in Section \ref{section:subcosets}.

\subsection*{Acknowledgments.} This work was supported by Harding Distinguished Postgraduate Scholars Programme. The author would like to thank Julia Wolf for many discussions and advice given in the course of this research, as well as Sean Prendiville, Julian Sahasrabudhe\cmnt{, and the anonymous reviewers} for their helpful comments on \cmnt{earlier versions} of this paper.

\section{Preliminaries}

This section sets out some essential definitions and concepts used in the rest of the paper, beginning with the formal definition of the density of a pattern.

\begin{definition}[Pattern density]
\label{def:pattern-density}
Given a system of linear forms $\mathcal{L} = \set{L_1, \ldots, L_m}$ and functions $f_1, \ldots, f_m:\vecsp{p}{n} \rightarrow [-1,1]$, define the operator
$$\Lambda_{\mathcal{L}}(f_1, \ldots, f_m) = \expct_{\underline{x} \in (\vecsp{p}{n})^{\cmnt{\ell}}} f_1(L_1(\vecline{x})) \ldots f_m(L_m(\vecline{x})).$$
For an $r$-colouring $\phi$, the \emph{density of $\mathcal{H} = (\mathcal{L}, \mathcal{X})$ under $\phi$} is given by
$$\Lambda_{\mathcal{H}}(\phi) = \sum_{\chi \in \mathcal{X}}\Lambda_{\mathcal{L}}(\id{\phi^{-1}(\chi(1))}, \ldots, \id{\phi^{-1}(\chi(m))}).$$
\end{definition}

\noindent It is easy to check that $\Lambda_{\mathcal{L}}$ satisfies the telescoping identity
\begin{equation}
    \label{eq:telescoping}
    \Lambda_{\mathcal{L}}(f_1, \ldots, f_m) - \Lambda_{\mathcal{L}}(g_1, \ldots, g_m) = \sum_{i=1}^m \Lambda_{\mathcal{L}}(h_1^\bracketed{i}, \ldots, h_m^\bracketed{i})
\end{equation}
where $h_j^\bracketed{i}$ is equal to $f_j$ when $j < i$, $f_i-g_i$ when $j=i$, and $g_j$ otherwise.

An integral tool in proving removal lemmas is an arithmetic analogue of Szemer\'edi's regularity lemma. Given a function $f:\vecsp{p}{n} \rightarrow [0,1]$, this arithmetic regularity lemma \cite{tower-type} provides a partition of $\vecsp{p}{n}$ into cosets of a `large' subspace such that $f$ behaves `pseudorandomly' on almost all cosets, with the measure of pseudorandomness given by Fourier uniformity. In the following definition, $\mathcal{P}(H)$ denotes the partition of $\vecsp{p}{n}$ into cosets of a subspace $H$.

\begin{definition}[Fourier uniformity]
\label{def:fourier}
    Let $H$ be a subspace of $\vecsp{p}{n}$. Given a function $F: \vecsp{p}{n} \rightarrow \mathbb{C}$ and elements $c, r \in \vecsp{p}{n}$, the \emph{Fourier transform of $F$ at $r$ on $H+c$} is defined as
        $\widehat{F{|_{H+c}}}(r) = \expct_{x \in H+c} F(x) e_p(r^T x),$
    where  $e_p(\cdot)$ denotes $\exp(2\pi i \cdot/p)$.
        
    A function $f: \vecsp{p}{n}\rightarrow \mathbb{C}$ is said to be \emph{$\epsilon$-uniform on $H+c$} if $\abs{\widehat{F{|_{H+c}}}(r)} \leq \epsilon$ for all $r \in \vecsp{p}{n}$, where $F = f - \expct_{x \in H+c} f(x)$. When $f$ is $\epsilon$-uniform on $H+c$, the latter is referred to as an \emph{$\epsilon$-regular coset for $f$}.

    The partition $\mathcal{P}(H)$ is \emph{$\epsilon$-regular for $f$} if for all but an $\epsilon$-proportion of $c \in \vecsp{p}{n}$, $H+c$ is $\epsilon$-regular for $f$.
\end{definition}

\noindent The precise statement of the arithmetic regularity lemm\cmnt{a}, due to Green \cite{tower-type}, can now be stated precisely, as follows.

\begin{theorem}[Arithmetic regularity lemma \cite{tower-type}]
\label{theorem:arl}
For all $\epsilon > 0$ and integer $r > 0$, there exists $C_{arl}(\epsilon, r)$ such that the following holds for any functions $f_1, \ldots, f_r:\vecsp{p}{n} \rightarrow [0,1]$. Given a subspace $H_0 \leqslant \vecsp{p}{n}$, there is a subspace $H \leqslant H_0$ of codimension at most $C_{arl}(\epsilon, r)$ in $H_0$ such that $\mathcal{P}(H)$ is $\epsilon$-regular for $f_1, \ldots, f_r$.
\end{theorem}

Fourier uniformity allows us to count the number of instances for linear systems of complexity 1. Specifically, the complexity of a linear system here refers to the true complexity of Gowers and Wolf \cite{gowers-wolf}, originally defined in terms of Gowers uniformity norms $\norm{f}_{U^{s+1}}$ \cite{gowers-norm}. While the latter will not be defined here, it is a well-known fact that functions with small $U^2$-norm are precisely those that have small Fourier transforms: in fact, it is easy to show that $\norm{f}_{U^2}^4 = \norm{\hat{f}}_4^4$ \cite[Lemma 2.4]{gowers-decompositions}, from which the last part of the definition below follows.

\begin{definition}[True complexity]
\label{def:complexity}
Let $\mathcal{L} = \set{L_1, \ldots, L_m}$ be a system of linear forms. The \emph{true complexity}  of $\mathcal{L}$ is the least positive integer $s$ (if it exists) with the following property. For every $\delta > 0$, there is an $\epsilon(\delta) > 0$ such that for any functions $F_1, \ldots, F_m: \vecsp{p}{n} \rightarrow [-1,1]$ satisfying $\min_i \norm{F_i}_{U^{s+1}} \leq \epsilon(\delta)$, $\babs{\Lambda_{\mathcal{L}}(F_1, \ldots, F_m)} \leq \delta.$

In particular, $\mathcal{L}$ has \emph{true complexity 1} if for every $\delta > 0$, there is an $\epsilon_{count}(\delta) > 0$ such that for any $F_1, \ldots, F_m: \vecsp{p}{n} \rightarrow [-1,1]$ with $\min_i \max_{\eta \in \vecsp{p}{n}} \abs{\widehat{F_i}(\eta)} \leq \epsilon_{count}(\delta)$,
$\babs{\Lambda_{\mathcal{L}}(F_1, \ldots, F_m)} \leq \delta.$
\end{definition}

\begin{example}
    The linear system $\set{x, x+d, x+2d}$, corresponding to 3-term arithmetic progressions, has true complexity 1 (for example, see \cite[Proposition 1.8]{montreal}). On the other hand, $\set{x, x+d, x+2d, x+3d}$ does not. There is a beautifully simple criterion that may be used to verify this fact. \cmnt{This criterion was conjectured and partially proved by Gowers and Wolf \cite[Theorem 6.1]{gowers-wolf-higher-degree} before being established in full generality by Hatami, Hatami, and Lovett \cite[Theorem 3.14]{hatami-lovett-complexity}. It states that} the true complexity of a linear system $\set{L_1, \ldots, L_m}$ over $\vecsp{p}{n}$ is equal to the smallest $s$ such that the set $\set{L_1^{s+1}, \ldots, L_m^{s+1}}$ is linearly independent. It is easy to check that while $x - 2(x+d) + (x+2d) = 0$, the set $\set{x^2, (x+d)^2, (x+2d)^2}$ is linearly independent, so the true complexity of the 3-term arithmetic progression is 1. On the other hand, $\set{x^2, (x+d)^2, (x+2d)^2, (x+3d)^2}$ satisfies $x^2 - 3(x+d)^2 + 3(x+2d)^2 - (x+3d)^2 = 0,$
    meaning that the true complexity of the 4-term arithmetic progression is at least 2.
\end{example}

Together with the telescoping identity \eqref{eq:telescoping} and an $\epsilon$-regular partition provided by Theorem \ref{theorem:arl}, Definition \ref{def:complexity} can be used to count instances of a linear system $\mathcal{L}$ that satisfy $L_i(\vecline{x}) \in H+c_i$ for given cosets $H+c_1, \ldots, H+c_m$.  In fact, if each $H+c_i$ is regular for $f_i$, then there are many such instances of $\mathcal{L}$, so long as the cosets $H+c_1, \ldots, H+c_m$ `align' in the right way. This is a necessary condition as, for example, there cannot be any instances of a $3$-term arithmetic progression $x_1, x_2, x_3$ with each $x_i \in H+c_i$ unless $H+c_1, H+c_2, H+c_3$ themselves form an arithmetic progression. The following definition arises as a special case of \cite[Definition 3.17]{full-induced}.

\begin{definition}[Coset consistency]
\label{def:consistency}
    Let $H \leqslant \vecsp{p}{n}$ be a subspace and let $\mathcal{L} = (L_1, \ldots, L_m)$ be a linear system of complexity 1 in $\ell$ variables. Given $c_1, \ldots, c_m \in \vecsp{p}{n}$, the cosets $H+c_1, \ldots, H+c_m$ are said to be \emph{consistent with $\mathcal{L}$} if there exists $\vecline{x} \in (\vecsp{p}{n})^{\cmnt{\ell}}$ such that $L_i(\underline{x}) \in H+c_i$ for each $i \in [m]$.
\end{definition}

If the cosets $H+c_1, \ldots, H+c_m$ are in fact consistent with $\mathcal{L}$, then the expected number of instances of $\mathcal{L}$ satisfying $L_i(\vecline{x}) \in H+c_i$ is approximately the same as one would expect in a random setting. This result, known as a counting lemma, appears in various forms in the literature, often stated in the more general setting of higher-order Fourier analysis (see, for instance, \cite[Theorem 3.10]{trans-invariant}, \cite[Theorem 3.19]{full-induced}, and the discussion following \cite[Lemma 1.3]{montreal}) or for linear systems of a particular form (such as \cite[Proposition 6.2]{tower-type}). The statement below pertains to general linear systems of complexity 1, with a proof given for completeness in Appendix \ref{appendix:counting-lemma}.

\begin{lemma}[Counting Lemma]
\label{lemma:counting-lemma}
Fix $\delta > 0$, $0 < \epsilon \leq \epsilon_{count}(\delta)$, and a linear system $\mathcal{L} = (L_1, \ldots, L_m)$ of complexity 1. Let $H \leqslant \vecsp{p}{n}$ be a subspace of codimension $d$. Then for any functions $f_1, \ldots, f_m:\vecsp{p}{n} \rightarrow [-1,1]$ and any $c_1, \ldots, c_m \in \vecsp{p}{n}$, the following holds. If the cosets $H+c_1, \ldots, H+c_m$ are consistent with $\mathcal{L}$ and are $\epsilon$-regular for $f_1, \ldots, f_m$, then
$$\Big|\Lambda_{\mathcal{L}}(f_1 \id{H+c_1}, \ldots, f_m \id{H+c_m}) - p^{-d(m-\rank(\mathcal{L}))}\prod_{i=1}^m \alpha_i\Big| \leq p^{-d(m-\rank(\mathcal{L}))} m \delta,$$
where $\alpha_i$ denotes the average of $f_i$ on $H+c_i$.
\end{lemma}

\section{Overview of the argument}
\label{section:argument-overview}
It is instructive to consider first the proof a non-induced arithmetic removal lemma. The latter states that if a given set doesn't have too many instances of a linear system $\mathcal{L}$, then all such instances can be eliminated by removing at most an $\epsilon$-proportion of $\vecsp{p}{n}$.

\begin{theorem}[Arithmetic removal lemma \cite{tower-type}]
\label{thm:arithmetic-removal}
    Fix $\epsilon > 0$ and a linear system $\mathcal{L}$ of complexity 1. There exists $\delta = \delta(\epsilon, \mathcal{L})$ satisfying the following.  If $A \subseteq \vecsp{p}{n}$ is a set such that $\Lambda_{\mathcal{L}}(\id{A}) \leq \delta$, then $A$ can be made $\mathcal{L}$-free on the whole of $\vecsp{p}{n}$ by removing at most $\epsilon \abs{\vecsp{p}{n}}$ elements from $A$.
\end{theorem}

\begin{figure}[h!]
     \centering
     \begin{subfigure}[t]{0.48\textwidth}
         \centering
         \includegraphics[alt={A square subdivided into nine equal subsquares, with seven blue and two grey. Some of the area of the square is hashed, with the boundary lines `wiggly' to represent an unspecified set.}, width=0.7\textwidth]{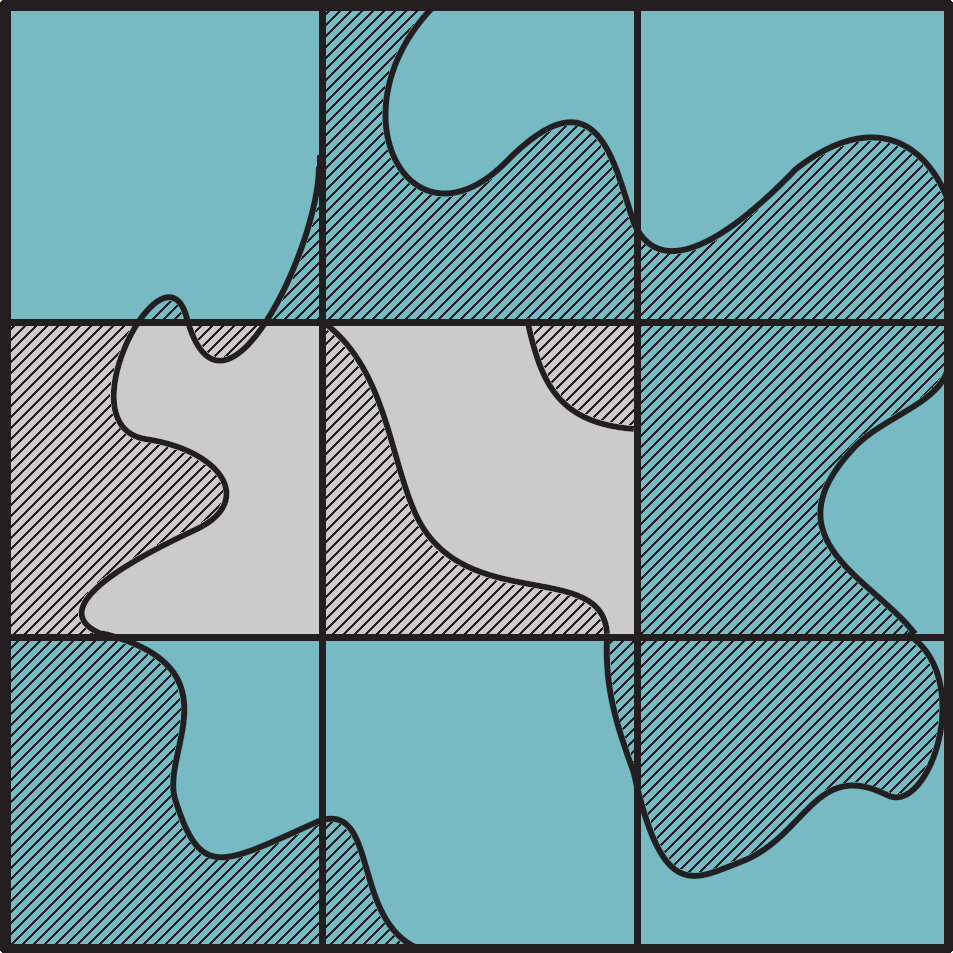}
     \end{subfigure}
     \hfill
     \begin{subfigure}[t]{0.48\textwidth}
         \centering
         \includegraphics[alt={An identically subdivided square, but this time all hashed areas inside the grey subsquares have been removed, as well as the hashed areas inside those blue subsquares where the intersection with the set was only small.}, width=0.7\textwidth]{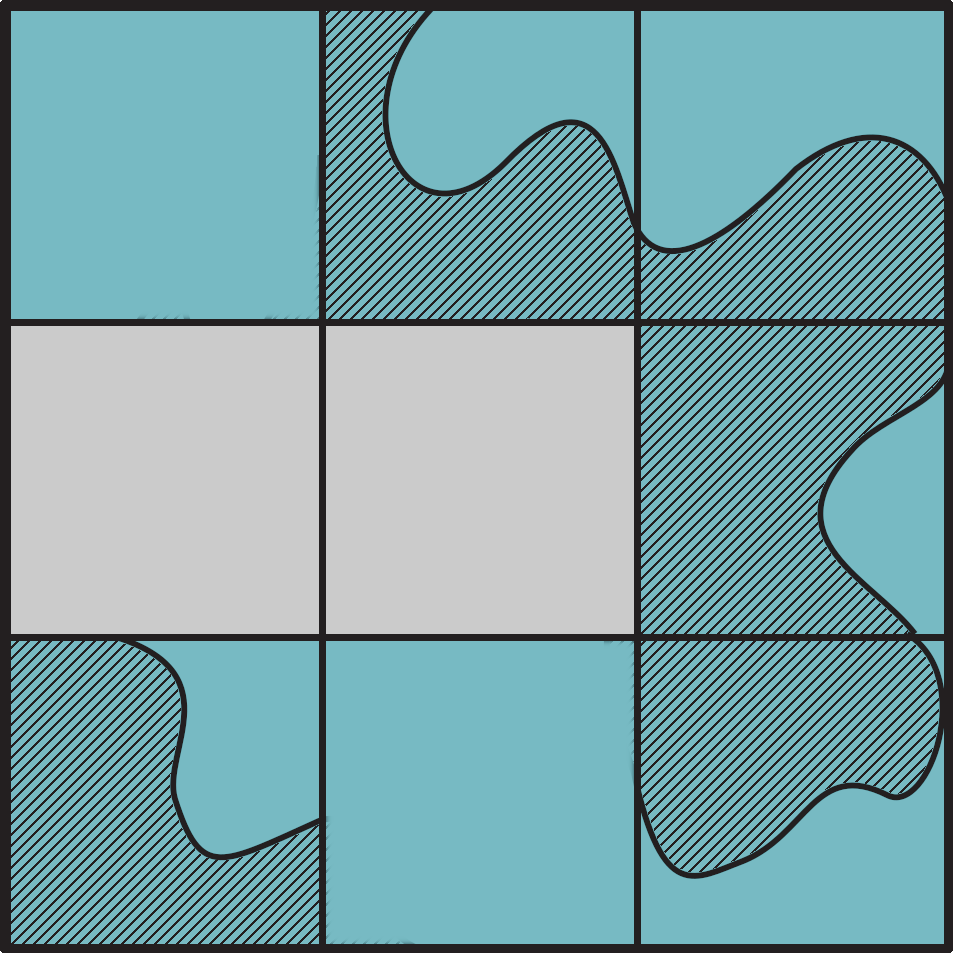}
     \end{subfigure}
     \caption{Sets $A$ (left) and $A'$ (right) depicted as hashed areas. Blue squares correspond to regular cosets.}
     \label{fig:removal-lemma}
\end{figure}

The proof of Theorem \ref{thm:arithmetic-removal} begins with an application of the arithmetic regularity lemma (Theorem \ref{theorem:arl}) to $\id{A}$, which gives a partition of $\vecsp{p}{n}$ into cosets of a subspace such that all but an $\epsilon$-proportion of the cosets are regular for $\id{A}$. The next step is to remove from $A$ all elements that lie either in a non-regular coset or in a coset on which $A$ has low density \cmnt{(see Figure \ref{fig:removal-lemma})}. Call the resulting set $A' \subseteq A$.

\begin{figure}[h]
     \centering
     \begin{subfigure}[t]{0.48\textwidth}
         \centering
         \includegraphics[alt={A square subdivided into nine subsquares, with the same colours and hashed area as the second square in Figure 1. There are three dots, each within the shaded area in a different subsquare. These dots are connected by edges, forming a triangle.}, width=0.7\textwidth]{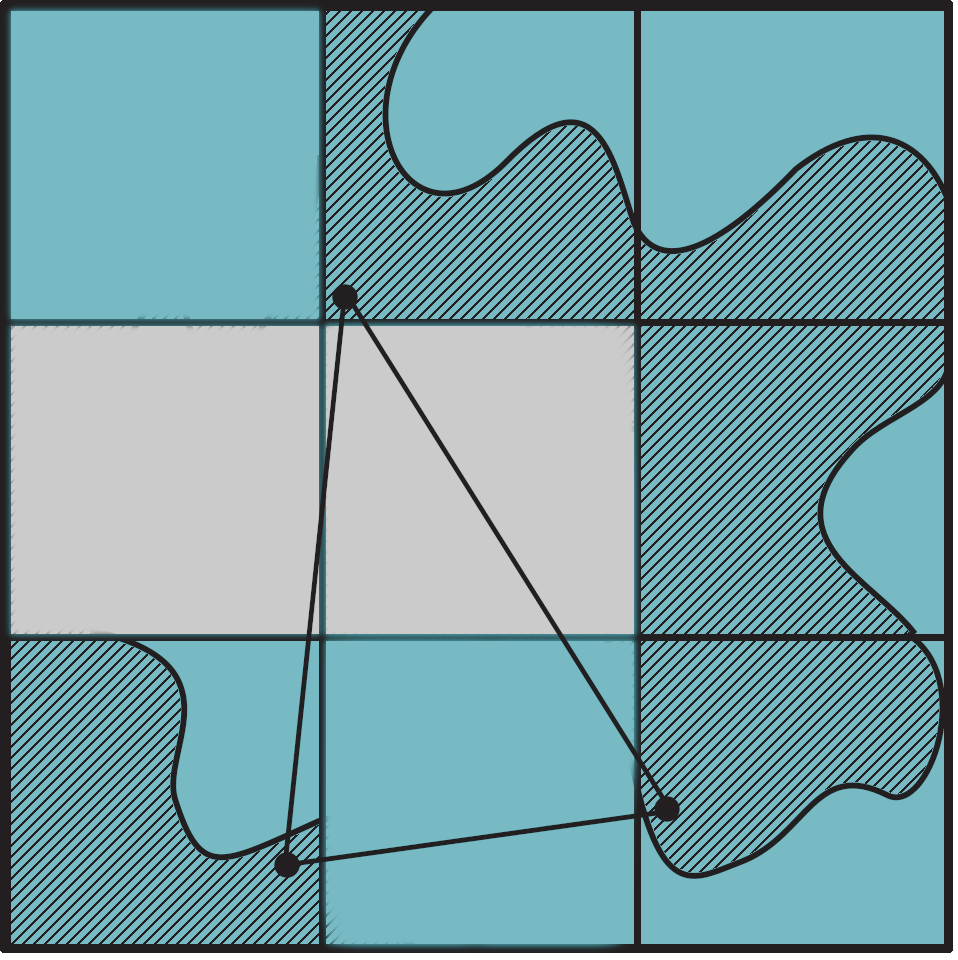}
     \end{subfigure}
     \hfill
     \begin{subfigure}[t]{0.48\textwidth}
         \centering
         \includegraphics[alt={A copy of the previous square but this time, there are many triangles, with vertices in the same subsquares.}, width=0.7\textwidth]{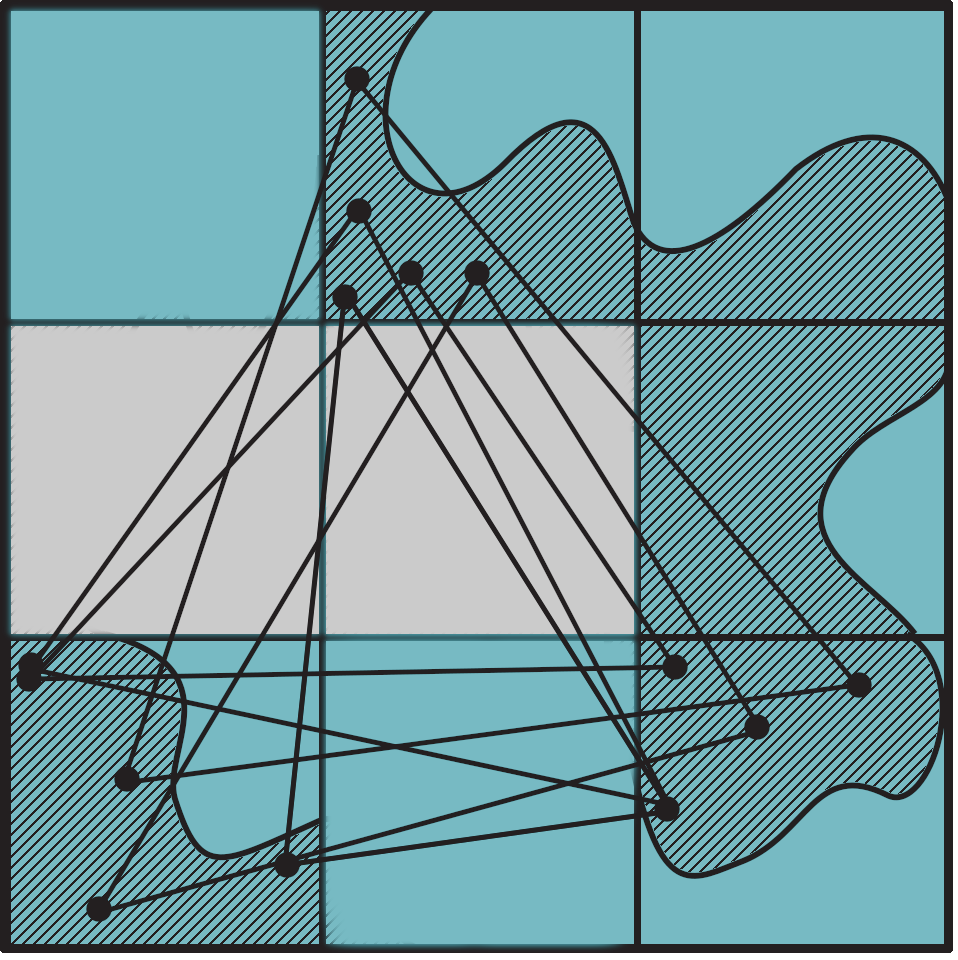}
     \end{subfigure}
     \caption{An instance of $\mathcal{L}$ in $A'$ would imply \textit{many} instances of $\mathcal{L}$ in $A$ as a consequence of Lemma \ref{lemma:counting-lemma}.}
     \label{fig:recolouring-contradiction}
\end{figure}

\noindent A crucial property of $A'$ is that if a coset contains at least one element of $A'$, then $A$ must both be Fourier-uniform and have high density on this coset. As a result, if there were a single instance of $\mathcal{L}$ in $A'$, we would be able to deduce that there are \textit{many} instances of $\mathcal{L}$ in $A$ via the counting lemma (Lemma \ref{lemma:counting-lemma}), as depicted in Figure \ref{fig:recolouring-contradiction}.

A set $A$ may equivalently be viewed as a 2-colouring of $\vecsp{p}{n}$ given by $\id{A}$, and the modified set $A'$ as an $\mathcal{L}$-free recolouring. In this way, Theorem \ref{thm:arithmetic-removal} is a special case of the induced removal lemma, where the pattern in question is $\mathcal{H} = (\mathcal{L}, \set{\chi})$ for the constant colouring $\chi \equiv 1$. Removing an element from $A$, then, corresponds to recolouring it with the colour $0$. Note that this never creates new instances of $\mathcal{H}$, which gives an easy way of avoiding non-regular cosets. On the other hand, if we had $\chi(1) = 0$ and $\chi(i) = 1$ everywhere else, the same would no longer hold.

This necessitates a different approach for proving induced removal lemmas. The existing proofs \cite{induced-over-f2, trans-invariant, induced-1, full-induced} all roughly follow the same strategy, involving two nested partitions rather than one.

\begin{enumerate}
    \item Use a strong version of the arithmetic regularity lemma (\cite[Theorem 5.4]{induced-1}) to find subspaces $H_2 \leqslant H_1$ such that
    \begin{enumerate}[label=(\roman*)]
        \item nearly all cosets of $H_2$ are regular for the given colouring $\phi$ (that is, for the indicator functions of all of its colour classes);
        
        \item the colour densities on almost all cosets of $H_2$ are a good approximation for the densities on the corresponding coset of $H_1$.
    \end{enumerate}
    
    \item Inside each coset of $H_1$, select a coset of $H_2$ to act as its representative in such a way that
    \begin{enumerate}[label=(\roman*)]
        \item each of the chosen subcosets is regular for $\phi$;
        \item the colour densities on (almost all) chosen subcosets approximate the densities on the cosets they are representing;
        \item the chosen subcosets preserve consistency with the given pattern, i.e.~if a collection of cosets of $H_1$ is consistent with the pattern, their subcoset representatives should be as well. 
    \end{enumerate}
\end{enumerate}

\noindent The resulting collection of subcosets can be thought of as a `regular model' for $H_1$ (see \cite[Section 2]{induced-1}): while it is impossible to guarantee that every coset of $H_1$ is regular for $\phi$ \cite{green-sanders}, we can pass to the chosen regular subcosets without losing too much information.

\begin{enumerate}
    \setcounter{enumi}{2}
    \item Define a recolouring $\psi$ as follows. For each coset of $H_1$, its subcoset representative determines whether a colour is `removed'. Specifically, if a colour occurs with low density on the subcoset, then it is replaced with some high-density colour on the whole of the coset. This recolours only a small proportion of all elements since the colour densities on the subcoset are approximately the same as on the coset it is representing.
\end{enumerate}

\noindent The resulting colouring $\psi$ has the property that if a colour occurs at least \textit{once} in a coset of $H_1$ under $\psi$, then it occurs with \textit{high density} on its subcoset representative under the original colouring $\phi$. Now, if there is an instance of the pattern under $\psi$, we can obtain a contradiction as in Figure \ref{fig:recolouring-contradiction} using the counting lemma on the subcoset representatives. This works since the representatives were chosen to preserve consistency (see Step 2).

\begin{figure}[h]
     \centering
     \begin{subfigure}[t]{0.48\textwidth}
         \centering
         \includegraphics[alt={A square, subdivided into nine intermediately-sized squares, each of which is subdivided into four further subsquares. There is a hashed area with `wiggly' boundaries to represent an unspecified set. In each intermediate square, one of the subsquares is outlines in orange. All but one of these are colour blue, and most of them contain a small black star.}, width=0.7\textwidth]{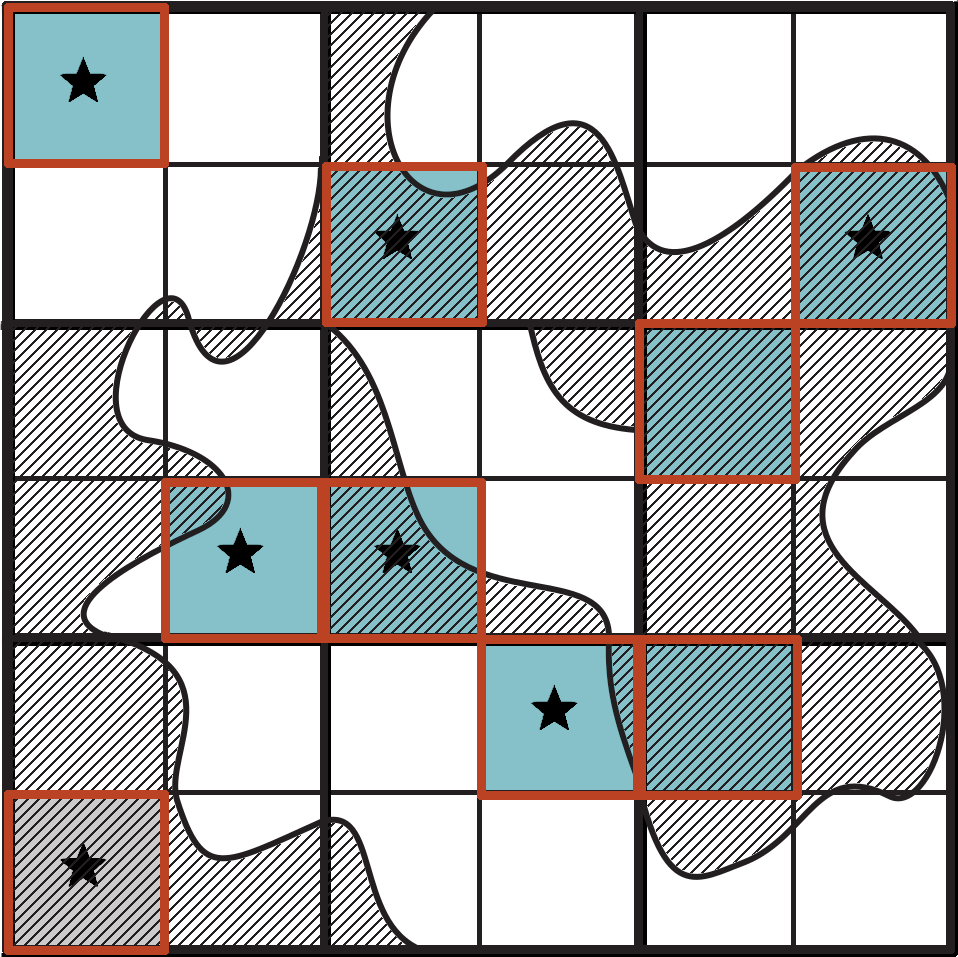}
         \caption{Subcoset representatives outlined in orange. Blue subcosets are regular for $\phi$, and the density is well-approximated on subcosets with a star.}
         \label{fig:selected}
     \end{subfigure}
     \hfill
     \begin{subfigure}[t]{0.48\textwidth}
         \centering
         \includegraphics[alt={An identical square but some of the hashed areas have been removed, namely: any areas in intermediate squares where the outlined subsquare did not contain a star; and any areas in intermediate squares whose intersection with the outlined subsquare was small. Additionally, the bottom left intermediate square has been made blank, with the words `Ramsey Patching Argument' written across it.}, width=0.7\textwidth]{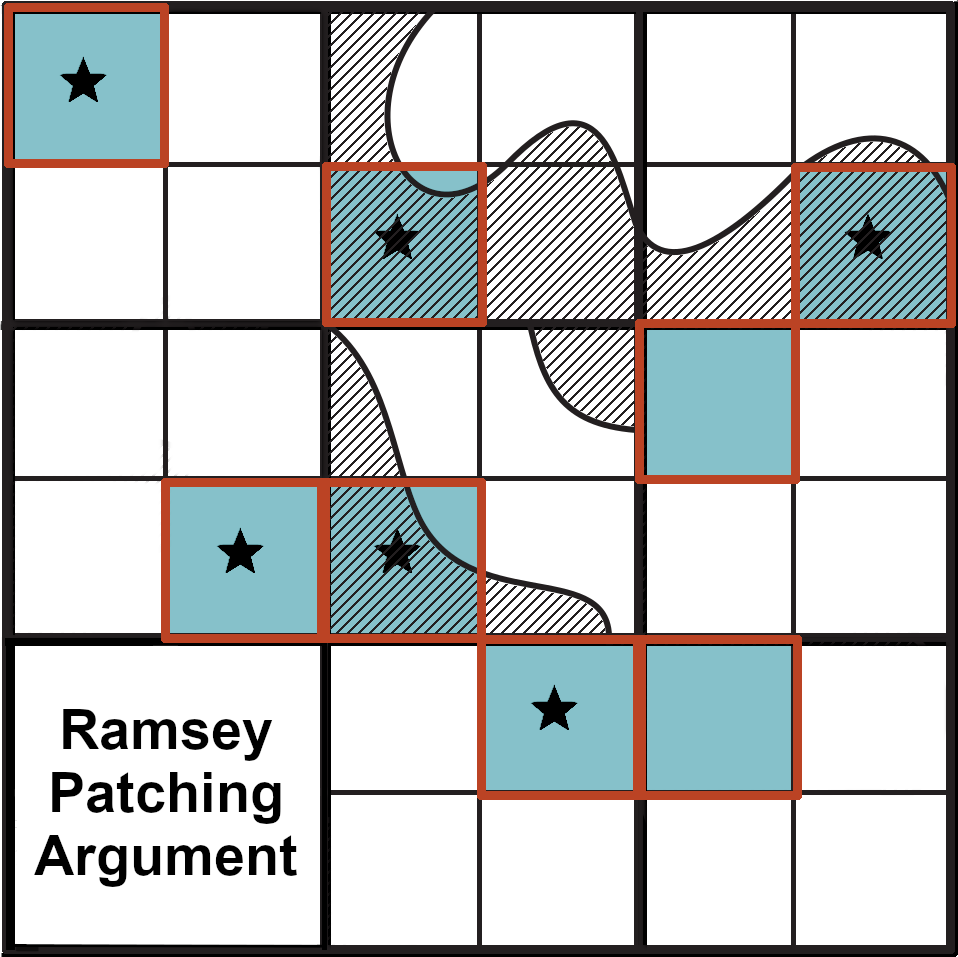}
         \caption{Each non-zero coset is recoloured according to its chosen subcoset; the zero coset is handled separately with a Ramsey-type argument.}
         \label{fig:recoulouring}
     \end{subfigure}
        \caption{The strategy of Fox, Tidor, and Zhao \cite{induced-1} for the induced arithmetic removal of complexity-1 patterns, depicted for a 2-colouring. Instances of the pattern which contain 0 may remain.}
        \label{fig:previous-approach}
\end{figure}

In order to make sure that the chosen subcosets preserve consistency, Fox, Tidor and Zhao \cite{induced-1} require that these subcosets themselves form a subspace. Of course, this leaves no choice for the subcoset representing the zero coset $H_1$, as it must be $H_2$ itself. This presents an issue since it is impossible to guarantee that $H_2$ is regular \cite{green-sanders}. As a result, their argument requires the zero coset to be treated separately from the rest of the space with a Ramsey-theoretic patching argument (Figure \ref{fig:recoulouring}). It is this last step that leads to a small number of pattern instances remaining after the recolouring in Theorem \ref{thm:induced-removal-general-1}.

This obstacle does not arise in the translation-invariant setting \cite{trans-invariant}, since the subspace of subcosets can be translated as desired without breaking consistency with the pattern. This added flexibility allows Bhattacharyya et al.~\cite{trans-invariant} to ensure that all chosen subcosets are regular.

\begin{figure}[h]
     \centering
     \begin{subfigure}[t]{0.45\textwidth}
         \centering
         \includegraphics[alt={Exactly the same square as in Figure 3(a).}, width=0.8\textwidth]{figure3a.png}
         \caption{A selection of cosets of $H_2$ inside those of $H_1$, for $H_2 \leqslant H_1$, as produced in \cite{induced-1}.}
         \label{fig:rado-intermediate}
     \end{subfigure}
     \hfill
     \begin{subfigure}[t]{0.45\textwidth}
         \centering
         \includegraphics[alt={The colouring, outlines and stars have been removed from the previously outlined subsquares. Instead, each subsquare has been subdivided even further into sixteen sub-subsquares. In each previously-outlined subsquare, there is now a triplet of sub-subsquares outlined in orange and all shaded blue. The choice of the triplet is exactly the same in each case. Six out of nine of these triplets have small black stars in each of the corresponding sub-subsquares, but three triplets do not.}, width=0.8\textwidth]{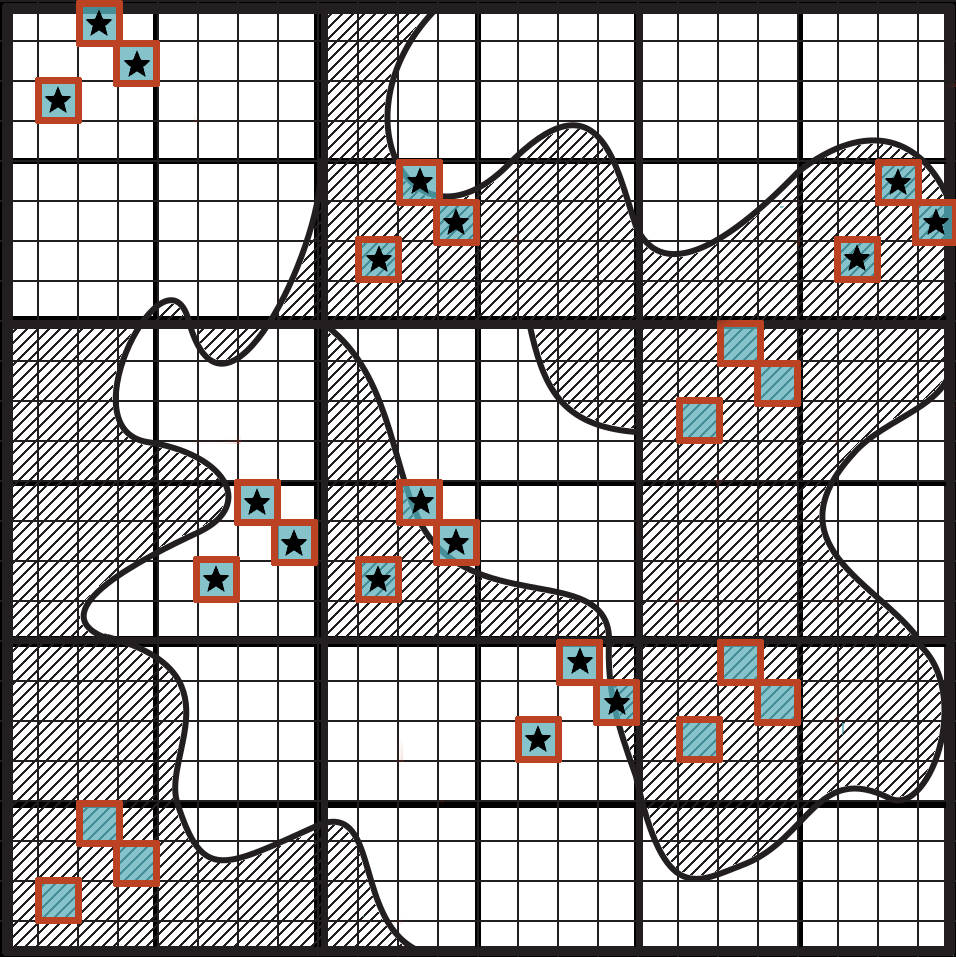}
         \caption{A selection of cosets of $H_3$ inside previously chosen cosets of $H_2$, for $H_3 \leqslant H_2$.}
         \label{fig:rado-selected}
     \end{subfigure}
     \caption{\cmnt{Subcoset selection scheme for partition-regular patterns}}
\end{figure}

While the approach outlined above does not apply directly to partition-regular patterns, some flexibility can still be regained by applying a modified subcoset selection strategy. The idea is to choose multiple subcosets inside each coset of $H_1$, so that the recolouring of each coset is governed by a group of representatives. \cmnt{An example of such a subcoset selection scheme is depicted in Figure \ref{fig:rado-selected}.}

\cmnt{For instance, consider an arithmetic pattern whose underlying linear system encodes solutions to $x+y=z$, which is partition-regular. We know there exists a subcoset selection as in Figure \ref{fig:selected}, so that there is some $H_2 \leqslant H_1$ and a subspace $U$ such that the chosen subcosets are given by $\set{H_2 + u: u \in U}$, and $(H_2+u_1) + (H_2+u_2) = H_2+u_3$ if and only if $(H_1+u_1)+(H_1+u_2)=H_1+u_3$. Now introduce a third regularity partition, given by a subspace $H_3 \leqslant H_2$. By defining an appropriate auxiliary colouring, we can use partition regularity to find cosets $W_x, W_y, W_z$ of $H_3$ in $H_2$ which all share some desirable properties while also satisfying $W_x + W_y = W_z$. In fact, we can find relatively many such triples by \cite[Theorem 2.1]{serra-vena}, which gives the required flexibility of choice. On the other hand, for each $u \in U$, $W_x+u, W_y+u, W_z+u$ lie in $H_2+u$, and as a result, $(W_x+u_1) + (W_y+u_2) = W_z+u_3$ if and only if $(H_1+u_1) + (H_1+u_2) = H_1+ u_3$, i.e.~consistency with the pattern is preserved.}

This \cmnt{approach} is reminiscent of the proof of the induced graph removal lemma by Alon, Fischer, Krivelevich and Szegedy (see the use of \cite[Corollary 3.4]{induced-graph} in the proof of \cite[Theorem 5.1]{induced-graph}). The technical details of \cmnt{the modified subcoset selection scheme and the resulting} proof of Theorem \ref{main} may be found in Section \ref{section:subcosets}.

\section{Linear systems and partition regularity}
\label{section:patterns}

Recall the well-known fact due to Rado \cite{rado} (extended to general abelian groups by Deuber \cite{deuber}) that a matrix $A$ is partition-regular if and only if $A$ satisfies Rado's column conditions.
\begin{definition}[Rado's column conditions]
\label{definition:column-condition}
Let $A$ be an $\cmnt{\ell \times} m$ matrix over $\vecsp{p}{}$. Then $A$ is said \cmnt{to} satisfy \emph{Rado's column conditions} if there are integers $0 < k_1 < \ldots < k_t=m$ and an ordering of the column vectors $\vecline{c}_1, \ldots, \vecline{c}_m$ such that
\begin{enumerate}[label=(\arabic*)]
    \item $\sum_{i=1}^{k_1} \vecline{c}_i = 0$;
    \item for each $1 < i \leq t$, $\sum_{j=k_{i-1}+1}^{k_i} \vecline{c}_j$ is in the span of $\vecline{c}_1, \ldots, \vecline{c}_{k_{i-1}}$.
\end{enumerate}
\end{definition}
As outlined in the previous section, the proof of Theorem \ref{main} relies on being able to choose subcosets with certain nice properties. Partition regularity facilitates this because we can define an auxiliary colouring on the set of all subcosets in a way that encodes these properties, and know that there are many monochromatic configurations of subcosets that are consistent with a given system $\mathcal{L}$ (see Lemma \ref{lemma:third-partition} for details). Specifically, we will be using the following quantitative form of Rado's theorem over finite fields due to Serra and Vena \cite[Theorem 2.1]{serra-vena}.

\begin{theorem}[Rado's theorem over finite fields \cite{serra-vena}]
\label{theorem:rado}
Fix integers $r, n > 0$, and let $A$ be an $\cmnt{\ell \times} m$ matrix over $\vecsp{p}{}$ that satisfies Rado's column conditions. Then there exist constants $c_{rado} = c_{rado}(p, r, m)$ and $n_{rado} = n_{rado}(p, r, m)$ such that for all $n \geq n_{rado}$ and for every $r$-colouring of $\vecsp{p}{n}$, the system $A\vecline{x}=0$ has at least $c_{rado}\abs{\vecsp{p}{n}}^{m-l}$ monochromatic solutions.
\end{theorem}

Of course, in order to apply Theorem \ref{theorem:rado} to arithmetic patterns, we need a way of turning a collection of linear forms into matrices. Such a translation is standard and outlined in this section for completeness.

Let $\mathcal{L} = \set{L_1, \ldots, L_m}$ be a system of linear forms in $\ell$ variables over $\vecsp{p}{}$ so that
$L_i(\underline{x}) = c_1^\bracketed{i} x_1 + \ldots + c_{\cmnt{\ell}}^\bracketed{i} x_{\cmnt{\ell}}.$
Define an $m \times \cmnt{\ell}$ matrix $\mathcal{M}(\mathcal{L}) = (c_j^\bracketed{i})$, where the $i$th row contains the coefficients of $L_i$. Then $\Img{\mathcal{M}}$ precisely corresponds to all values that the linear system $\mathcal{L}$ can take.

Moreover, it is easy to see that the linear dependencies between $L_1, \ldots, L_m$ are given by vectors $\underline{y} \in \vecsp{p}{m}$ satisfying $\underline{y} \mathcal{M}=0$, or, equivalently, $\underline{y} \in \ker(\mathcal{M}^T)$. On the other hand, these dependencies are precisely the vectors in $\Img{\mathcal{M}}^\perp$, as the following lemma shows.

\begin{lemma}
\label{lemma:matrices}
Let $A$ be any matrix. Then $\Img{A}^\perp=\ker(A^T)$.
\end{lemma}
\begin{proof}
Take $v \in \Img{A}^\perp$. Then for all $x$, $v^TAx=0$ and so $v^T A = 0$. Hence $A^T v = 0$, i.e.~$v \in \ker(A^T)$. For the other direction, take $v \in \ker(A^T)$. Then for all $x \in \vecsp{p}{n}$, $v^TAx=(A^Tv)^Tx=0.$ Hence $v \in \Img{A}^\perp$.
\end{proof}

\noindent Now let $\set{v_1, \ldots, v_{\cmnt{k}}}$ be a basis of $\ker(\mathcal{M}^T)$ and define $\mathcal{K}(\mathcal{L})$ to be the matrix with rows given by $\set{v_1^T, \ldots, v_{\cmnt{k}}^T}$. While the definition of $\mathcal{K}(\mathcal{L})$ depends on the choice of $\set{v_1, \ldots, v_{\cmnt{k}}}$, the properties we are interested in are invariant under a change of basis, so this choice is of no import.

\begin{prop}
\label{prop:solution-matrix}
The solutions to $\mathcal{K}(\mathcal{L}) \underline{z} = 0$ are precisely the values taken by $\mathcal{L}$ as $\underline{x}$ ranges over $\vecsp{p}{\cmnt{\ell}}$. That is, for each such $\underline{z}$ there is an $\underline{x} \in \vecsp{p}{\cmnt{\ell}}$ such that $\underline{z} = (L_1(\underline{x}), \ldots, L_m(\underline{x}))$.
\end{prop}
\begin{proof}
The definition of $\mathcal{K}(\mathcal{L})$ implies that $\Img{\mathcal{K}^T}= \ker(\mathcal{M}^T)$. Then, by Lemma \ref{lemma:matrices}, $\ker(\mathcal{K})=\Img{\mathcal{K}^T}^\perp= \ker(\mathcal{M}^T)^\perp=(\Img{\mathcal{M}}^\perp)^\perp = \Img{\mathcal{M}},$
which implies the result.
\end{proof}

\noindent In line with Proposition \ref{prop:solution-matrix}, we will say that a linear system $\mathcal{L}$ is partition-regular if and only if $\mathcal{K}(\mathcal{L})$ satisfies Rado's column conditions.

\section{Proof of Theorem \ref{main}}
\label{section:subcosets}
This section contains the core of the proof and the main original contribution of the paper. Whereas previous methods \cite{induced-1}, \cite{trans-invariant} depended on choosing a single suitable subcoset within each coset of the coarser partition (Figure \ref{fig:selected}), the proof of Theorem \ref{main} requires an additional refining step which allows us to choose \textit{several} subcosets at once. In fact, we will use the subcoset selection of \cite{induced-1} (Figure \ref{fig:rado-intermediate}) as an intermediate stage, then find the final subcosets inside each of the intermediate ones (Figure \ref{fig:rado-selected}). Here regular cosets are coloured in blue, while the star inside a coset is used to depict property \ref{prop:close-density} of Proposition \ref{prop:two-partitions} (which is a restatement of \cite[Proposition 3.2]{induced-1}).

\begin{prop}[Intermediate subcoset selection \cmnt{\cite{induced-1}}]
\label{prop:two-partitions}
Fix $\epsilon, \zeta > 0$ and integers $r, d > 0$. There exist $n_{reg} = n(\epsilon, \zeta, r, d)$ and $C_{reg}(\epsilon, \zeta, r, d)$ such that for any $n > n_{reg}$ the following holds. Given a subspace $H_0 \leqslant \vecsp{p}{n}$ of codimension $d$ and functions $f_1, \ldots, f_r: \vecsp{p}{n} \rightarrow [0,1]$, there are subspaces $H_2 \leqslant H_1$ of $H_0$ with codimensions $\codim(H_i) = D_i$ in $\vecsp{p}{n}$, as well as a choice of complement $U$ satisfying $U \oplus H_1 = \vecsp{p}{n}$ such that
\begin{enumerate}[label=(\roman*)]
    \item $D_1 \leq D_2 \leq C_{reg}(\epsilon, \zeta, r, d)$;
    \item for all $u \in U\backslash\set{0}$, $H_2+u$ is $\epsilon$-regular for $f_1, \ldots, f_r$;
    \item \label{prop:close-density} for all but at most a $\zeta$-proportion of $u \in U$ and all $i \in [r]$,
    $$\cmnt{\left|\expct_{x \in H_2+u}f_i(x) - \expct_{x \in H_1+u}f_i(x)\right| < \zeta.}$$
\end{enumerate}
\end{prop}

\noindent It is not hard to show that the property of being $\epsilon$-regular is inherited by subcosets and, moreover, provides us with a way of inheriting property \ref{prop:close-density} as well.

\begin{prop}[Regularity properties]
\label{prop:uniformity}
Fix $\epsilon > 0$. Let $f: \vecsp{p}{n} \rightarrow [0,1]$, and let $H_2 \leqslant H_1$ be subspaces of $\vecsp{p}{n}$ with $\mathrm{codim}_{H_1} H_2 = d$. Suppose that a coset $C_1$ of $H_1$ is $\epsilon$-regular for $f$. Then for any coset $C_2$ of $H_2$ such that $C_2 \subseteq C_1$,
\begin{enumerate}[label = (\roman*)]
    \item \label{prop:uniformityA} $C_2$ is $p^{d}\epsilon$-regular for $f$;
    \item \label{prop:uniformityB} $\babs{\expct_{x \in C_2} f(x) - \expct_{x \in C_1} f(x)} \leq p^{d}\epsilon$.
\end{enumerate}
\end{prop}

\begin{proof}
    Write $\alpha_i = \expct_{x \in C_i} f(x)$ and $C_i = H_i+z_i$ for $i=1,2$, where $z_1, z_2 \in \vecsp{p}{n}$. Since $C_1$ is $\epsilon$-regular for $f$, we know that for any $r \in \vecsp{p}{n}$,
    \begin{equation}
        \label{eq:f-is-regular}
        \babs{\expct_{x \in C_1} (f(x) - \alpha_1) e_p(r^T x)} \leq \epsilon.
    \end{equation}
    As $\mathrm{codim}_{H_1} H_2 = d$, there exist some $h_1, \ldots, h_d \in H_1$ such that $H_2 = H_1 \cap \gen{h_1, \ldots, h_d}^\perp$. The indicator function $\id{H_2}$ can then be written as
    \begin{equation}
        \label{eq:indicator}
        \id{H_2}(x) = \id{H_1}(x) \prod_{i=1}^d \expct_{a_i \in \vecsp{p}{}} e_p(a_i h_i^T x) = \id{H_1}(x) \expct_{a \in \vecsp{p}{d}} e_p(R_a^T x),
    \end{equation}
    where $R_a = \sum_{i=1}^d a_i h_i$. Now fix some $r \in \vecsp{p}{n}$. If $r \in H_2^\perp$, then $\widehat{f|_{C_2}}(r) = 0$ trivially. Otherwise,
    \begin{align*}
        \widehat{f|_{C_2}}(r) &= \frac{\abs{C_1}}{\abs{C_2}}\expct_{x \in C_1} (f(x)-\alpha_2) \id{H_2}(x-z_2) e_p(r^T x)\\
        &= p^d \expct_{a \in \vecsp{p}{d}} e_p(-R_a^T z_2) \Big[\expct_{x \in C_1} (f(x) - \alpha_2) e_p((r+R_a)^T x)\Big].
    \end{align*}
    Then using \eqref{eq:f-is-regular} and the triangle inequality gives
    $$\babs{\widehat{f|_{C_2}}(r)} \leq p^d \epsilon + \babs{\expct_{x \in C_1} (\alpha_1-\alpha_2) e_p((r+R_a)^T x)}.$$
    But the last term in this equation is 0 since $(r+R_a)^Tx \not\equiv 0$ on $C_1$ when $r \notin H_2^\perp$. This is because $H_1^\perp \leqslant H_2^\perp$ so $r \notin H_1^\perp$ and $R_a \in H_1$ by definition, which implies $r+R_a \notin H_1^\perp$. Hence part \ref{prop:uniformityA} of the proposition holds. Part \ref{prop:uniformityB} is proved similarly. From equation \eqref{eq:indicator},
    $$\expct_{x \in C_2} f(x) - \alpha_1 = p^d \expct_{a \in \vecsp{p}{d}} e_p(-R_a^T z_2) \Big[\expct_{x \in C_1} (f(x) - \alpha_1) e_p(R_a^T x)\Big],$$
    which implies, by \eqref{eq:f-is-regular} and the triangle inequality, that $\abs{\alpha_2 - \alpha_1} \leq p^d \epsilon$, as required.
\end{proof}

Thanks to Proposition \ref{prop:uniformity}, if $\codim_{H_2} H_3$ is sufficiently small, any cosets of $H_3$ picked inside the regular cosets of $H_2$ (as in Figure \ref{fig:rado-selected}) will still be regular, and $f$ will still have approximately the same density on them as on the coset of $H_1$ containing them. Since non-zero cosets of $H_2$ can be made regular by Proposition \ref{prop:two-partitions}, it is enough to make a `good' selection of subcosets inside the zero coset of $H_2$. This is the purpose of the following lemma.

\begin{lemma}[Selecting multiple subcosets]
\label{lemma:third-partition}
Fix $\epsilon, \delta > 0$, an integer $r > 0$, and an $\cmnt{\ell \times} m$ partition-regular matrix $A$ over $\vecsp{p}{}$. There exist $n_{zreg} = n(r, m)$ and $C_{zreg}(\epsilon, r, m)$ such that for any $n > n_{zreg}$ the following holds. Given functions $f_1, \ldots, f_r: \vecsp{p}{n} \rightarrow [0,1]$, there is a subspace $H \leqslant \vecsp{p}{n}$ of codimension $D \leq C_{zreg}$, as well as $z_1, \ldots, z_m \in \vecsp{p}{n}$ such that
    \begin{enumerate}[label=(\roman*)]
        \item  $A (z_1, \ldots, z_m)^T = 0$;
        \item \label{prop:third-partition-regular} for all $j \in [m]$, $H+z_j$ is $\epsilon$-regular for $f_1, \ldots, f_r$;
        \item \label{prop:third-partition-densities} for each $i \in [r]$, $\expct_{x \in H + z_j} f_i(x) < \delta$ for all $j \in [m]$ or $\expct_{x \in H + z_j} f_i(x) \geq \delta$ for all $j \in [m]$.
    \end{enumerate}
\end{lemma}

\begin{proof}
    \cmnt{Let $n_{zreg} = n_{rado}(p, 2^r, m)$, $\epsilon' = \min(\epsilon, c_{rado}(p, 2^r, m)/2)$, and $C_{zreg} = C_{arl}(\epsilon', r)+2n_{rado}(p, 2^r, m)$, where $n_{rado}$ and $c_{rado}$ are as in Theorem \ref{theorem:rado}, and $C_{arl}$ as in the arithmetic regularity lemma (Theorem \ref{theorem:arl}). Choose any subspace $H_0$ of codimension $n_{rado}(p, 2^r, m)$ and apply the arithmetic regularity lemma} to $f_1, \ldots, f_r$ with \cmnt{parameters} $\epsilon'$\cmnt{, $r$, and $H_0$} to obtain a subspace $H \leqslant H_0$ such that
    \begin{itemize}
        \item the codimension $D$ of $H$ in $\vecsp{p}{n}$ satisfies $n_{rado}(p, 2^r, m) \leq D \leq \cmnt{C_{zreg}}$;

        \item all but an $\epsilon'$-proportion of the cosets of $H$ are $\epsilon$-regular for $f_1, \ldots, f_r$.
    \end{itemize}

    Let $U$ be \cmnt{a subspace} such that $U \oplus H = \vecsp{p}{n}$, and define a colouring $\psi: U \rightarrow \mathcal{P}([r])$ by setting $\psi(u) = \set{i \in [r] \st \expct_{x \in H+u} f_i(x) \geq \delta}$, so that $\psi(u)$ corresponds to the set of those $f_i$ that have high density on $H+u$. Apply Theorem \ref{theorem:rado} to $\psi$ in $U$ (which is isomorphic to $\vecsp{p}{D}$) to find at least $c_{rado}(p, 2^r, m) \abs{\vecsp{p}{D}}^{m-l}$ monochromatic solutions to $A \vecline{z} = 0$ under $\psi$. Note that by definition of $\psi$, any such solution satisfies property \ref{prop:third-partition-densities} of the lemma.

    \begin{figure}[h]
     \centering
     \includegraphics[alt={A square, subdivided into sixteen subsquares, all but two of them shaded blue. Each subsquare contains a label. From left to right in each of the four rows, these are: (b), (r,b,g), (r,g), (r,g); (b), (), (), (r,g,); (r,g), (r,g), (r,b,g), (g); (), (r,b,g), (r), (b).}, width=0.3\textwidth]{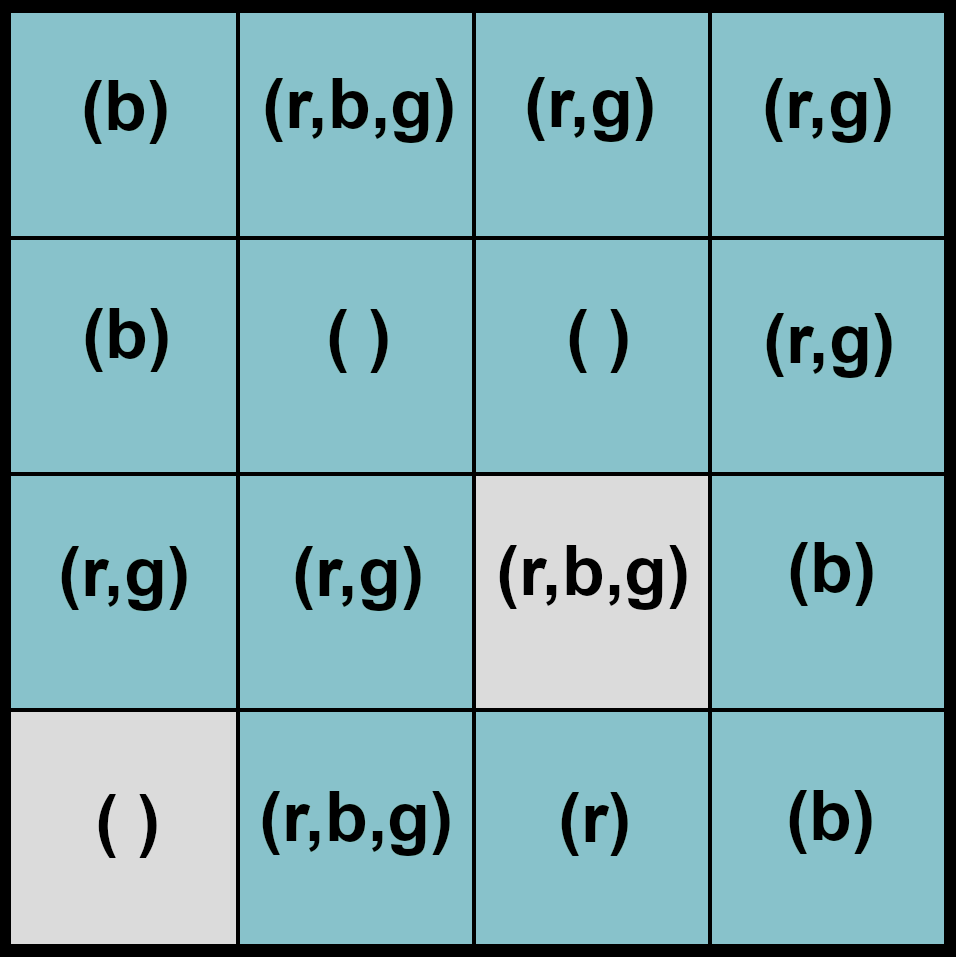}
     \hspace{0.15\textwidth}
     \includegraphics[alt={The same square but now three of the blue subsquares are outlined in orange -- each of them has the same label, namely (r,g).}, width=0.3\textwidth]{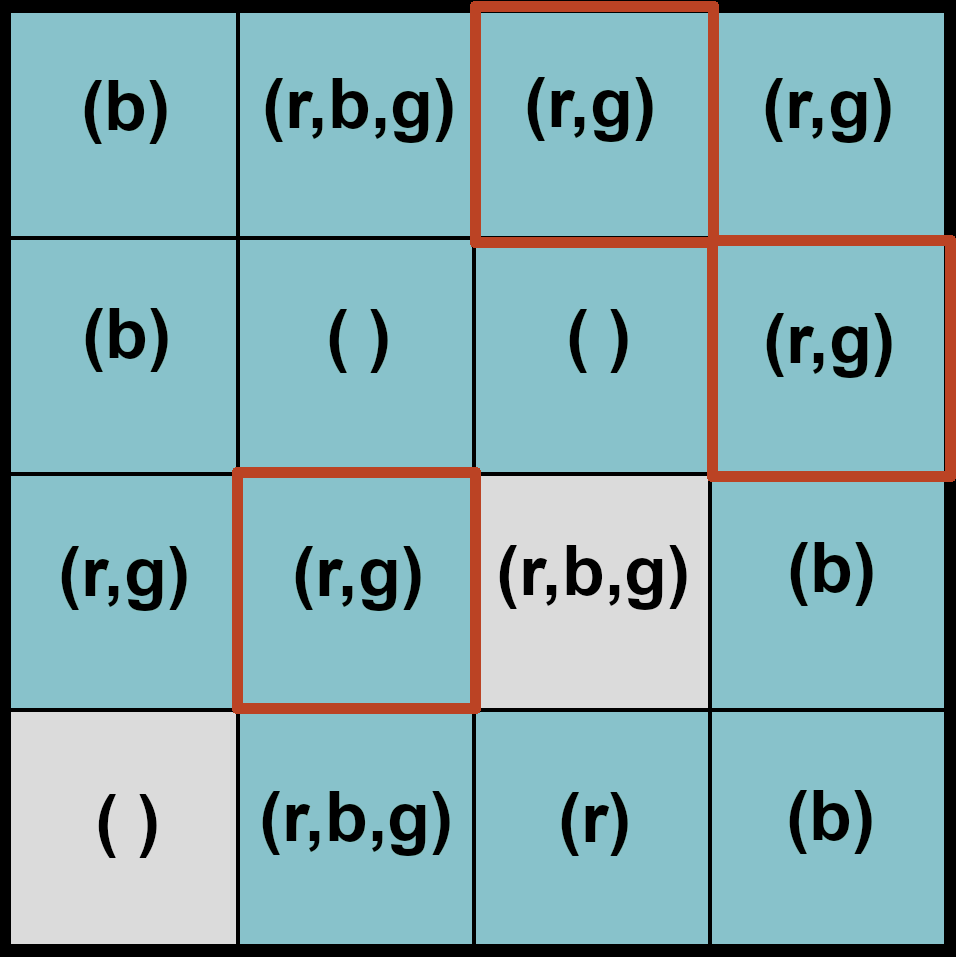}
        \caption{An auxiliary colouring $\psi$ arising in the proof of Lemma \ref{lemma:third-partition} from a 3-colouring $\phi:\vecsp{p}{n} \rightarrow \set{r,b,g}$. Due to a large proportion of regular cosets (shown in blue), we can find a monochromatic solution consisting entirely of regular cosets.}
        \label{fig:third-partition}
\end{figure}

    On the other hand, the number of solutions for which property \ref{prop:third-partition-regular} does not hold is at most $\epsilon'\abs{\vecsp{p}{D}}^{m-l} \leq \frac12 c_{rado}(p, 2^r, m)\abs{\vecsp{p}{D}}^{m-l}$, which is strictly less than the total number of monochromatic solutions. As such, there is at least one solution $(z_1, \ldots, z_m)$ that satisfies all properties\cmnt{, as depicted in Figure \ref{fig:third-partition}.}
\end{proof}

All results stated so far have referred to 1-bounded functions whereas our aim is to say something about an $r$-colouring $\phi:\vecsp{p}{n} \rightarrow [r]$. Of course, it is easy to turn $\phi$ into a collection of $r$ functions by letting $\phi_i:\vecsp{p}{n} \rightarrow [0,1]$ be the indicator function of colour $i$, i.e.~$\phi_i = \id{\phi^{-1}(i)}$. We can then combine the other results in this section to define an appropriate recolouring for Theorem \ref{main}.

\begin{lemma}[Recolouring]
\label{lemma:recolouring1}
Fix $\epsilon > 0$, an integer $r > 0$, and an $\cmnt{\ell \times} m$ partition-regular matrix $A$ over $\vecsp{p}{}$. There exist $n_{rcl} = n_{rcl}(\epsilon, r, m)$ and $C_{rcl} = C_{rcl}(\epsilon, r, m)$ such that for any $n > n_{rcl}$ the following holds. Given an $r$-colouring $\phi: \vecsp{p}{n} \rightarrow [r]$, there exist subspaces $H_3 \leqslant H_2 \leqslant H_1$ of $\vecsp{p}{n}$ with codimensions $\codim(H_i) = D_i$, a choice of complement $U$ satisfying $U \oplus H_1 = \vecsp{p}{n}$, and $z_1, \ldots, z_m \in H_2$ such that
    \begin{enumerate}[label=(\roman*)]
        \item $\log_p(\frac{\epsilon}{4r}) \leq D_1, D_2, D_3 \leq C_{rcl}$;
        \item $A (z_1, \ldots, z_m)^T = 0$;
        \item \label{prop:recolouring-regularity} for all $u \in U$ and $j \in [m]$, $H_3+u+z_j$ is $\epsilon$-regular for $\phi_1, \ldots, \phi_r$.
    \end{enumerate}
    
    Additionally there is an $r$-colouring $\phi': \vecsp{p}{n} \rightarrow [r]$ such that $\phi'$ differs from $\phi$ in at most $\epsilon \abs{\vecsp{p}{n}}$ places and satisfies the following property.
    \begin{enumerate}[label=(\roman*)]
        \setcounter{enumi}{3}
        \item \label{prop:recolouring-density} For any $u \in U$, if $\phi'(y) = i$ for some $y \in H_1+u$, then $\expct_{x \in H_3+u+z_j} \phi_i(x) \geq \epsilon/8r$ for all $j \in [m]$.
    \end{enumerate}
\end{lemma}

\begin{proof}
    \cmnt{With $n_{zreg}, C_{zreg}$ as in Lemma \ref{lemma:third-partition} and $n_{reg}, C_{reg}$ as in  Proposition \ref{prop:two-partitions},} define parameters $D_0 = \ceil{\log_p({\epsilon}/{4r})}$, $C_{zreg} = C_{zreg}\cmnt{(\epsilon, r, m)}$, $\epsilon' = \epsilon p^{-C_{zreg}}/{8r}$ and $\zeta = {\epsilon}/{4r}$; \cmnt{then $n_{rcl}$ and $C_{rcl}$ may be taken to be} $C_{rcl} = C_{reg}(\epsilon', \zeta, \cmnt{r}, D_0) + C_{zreg}$ and
    $$n_{rcl} = \max\left(
        n_{reg}(\epsilon', \zeta, r, D_0),
        n_{zreg}(r, m) + C_{reg}(\epsilon', \zeta, \cmnt{r}, D_0)
    \right).$$
    Let $H_0$ be an arbitrary subspace of codimension $D_0$, and apply Proposition \ref{prop:two-partitions} to $\phi_1, \ldots, \phi_r$ and $H_0$ with parameters $\epsilon'$ and $\zeta$. The result is subspaces $H_2 \leqslant H_1$ of $H_0$ with codimensions $\codim(H_i) = D_i$, as well as a complement $U$ satisfying $U \oplus H_1 = \vecsp{p}{n}$ such that the following holds:
    \begin{enumerate}[label=(\alph*)]
        \item $\log_p({\epsilon}/{4r}) \leq D_1, D_2 \leq C_{reg}(\epsilon', \zeta, \cmnt{r}, D_0)$;
        \item \label{prop:inter-regular} for all $u \in U\backslash\set{0}$, $H_2+u$ is $\epsilon'$-regular for $\phi_1, \ldots, \phi_r$;
        \item for all but at most a $\zeta$-proportion of $u \in U$ and all $i \in [r]$,
    $$\babs{\expct_{x \in H_2+u}\phi_i(x) - \expct_{x \in H_1+u}\phi_i(x)} < \zeta.$$
    \end{enumerate}

    Now apply Lemma \ref{lemma:third-partition} to $H_2$ with parameters $\epsilon$ and $\delta = \epsilon/8r$ to obtain a subspace $H_3 \leqslant H_2$ of codimension $D_3 = C_{zreg} + D_2 \leq C_{rcl}$ in $\vecsp{p}{n}$, and $z_1, \ldots, z_m \in H_2$ such that
    \begin{enumerate}[label=(\alph*)]
        \setcounter{enumi}{3}
        \item $A (z_1, \ldots, z_m)^T = 0$;
        \item \label{prop:zero-regular} for all $j \in [m]$, $H_3+z_j$ is $\epsilon$-regular for $\phi_1, \ldots, \phi_r$;
        \item \label{prop:same-colour-density} for all $i \in [r]$, either all $j \in [m]$ satisfy $\expct_{x \in H_3 + z_j} \phi_i(x) < \delta$ or all $j \in [m]$ satisfy $\expct_{x \in H_3 + z_j} \phi_i(x) \geq \delta$.
    \end{enumerate}

    Property \ref{prop:inter-regular} and Proposition \ref{prop:uniformity} imply that for all $u \in U \backslash \set{0}$ and $j \in [m]$, the coset $H_3+u+z_j$ is $\epsilon' p^{C_{zreg}}$-regular for $\phi_1, \ldots, \phi_r$. By the choice of parameters, $\epsilon' p^{C_{zreg}} \leq {\epsilon}/{8r}$ which together with \ref{prop:zero-regular} gives part \ref{prop:recolouring-regularity} of the lemma. All that remains is to define a suitable recolouring $\phi'$.

    Fix $u \in U$. If $u \neq 0$, let $c$ be a colour with density at least $1/r$ on $H_2+u$.
    \begin{itemize}
        \item If $\babs{\expct_{x \in H_2+u} \phi_i(x) - \expct_{x \in H_1+u} \phi_i(x)} > \zeta$ for some $i \in [r]$, recolour all of $H_1+u$ with $c$; note that this happens for at most a $\zeta$-proportion of all possible $u \in U \backslash \set{0}$.
        \item Otherwise, for each colour $i \in [r]$ with density less than $\epsilon/4r$ on $H_2+u$, recolour all occurrences of $i$ in $H_1+u$ with $c$; note that this changes at most a ${\epsilon}/{4r} + \zeta \leq {\epsilon}/{2r}$ proportion of $H_1+u$.
    \end{itemize}
    These steps affect at most a $\zeta + r {\epsilon}/{2r} \leq 3\epsilon/4$ proportion of the whole space.

    For $u = 0$, let $c$ be a colour that has density at least $1/r$ on $H_3+z_1$. Then by property \ref{prop:same-colour-density}, $c$ has density at least $\delta = \epsilon/8r$ on each $H_3+z_j$. Recolour all of $H_1$ with $c$ and note that this means $\phi'$ satisfies property \ref{prop:recolouring-density} when $u = 0$.

    In total, the recolouring $\phi'$ defined by following the steps above differs from $\phi$ in at most $(3\epsilon/4 + {p^{-D_1}})\abs{\vecsp{p}{n}} \leq \epsilon \abs{\vecsp{p}{n}}$ places.

    To show that property \ref{prop:recolouring-density} holds for all $u$, fix $u \in U \backslash \set{0}$ and some $y \in H_1+u$. By definition of $\phi'$, if $\phi'(y) = i$, then $\expct_{x \in H_2+u} \phi_i(x) \geq \epsilon/4r$. On the other hand, $H_3+u+z_j \subseteq H_2+u$ for any $j \in [m]$, and $H_2+u$ is $\epsilon'$-uniform for $\phi_i$; then Proposition \ref{prop:uniformity} implies
    $$\expct_{x \in H_3+u+z_j} \phi_i(x) \geq \expct_{x \in H_2+u} \phi_i(x) - \epsilon' p^{C_{zreg}} \geq \frac{\epsilon}{4r} - \frac{\epsilon}{8r} = \frac{\epsilon}{8r},$$
    as required. Together with the earlier observation that property \ref{prop:recolouring-density} holds for $u = 0$, this completes the proof.
\end{proof}

Lemma \ref{lemma:recolouring1} contains all of the ingredients necessary to prove Theorem \ref{main}, restated as Theorem \ref{theorem:main-linear} below. To complete the proof, we must show that the recolouring provided by Lemma \ref{lemma:recolouring1} has no instances of the given partition-regular pattern.

\begin{theorem}
    \label{theorem:main-linear}
    Fix $\epsilon > 0$, and let $\mathcal{H} = (\mathcal{L}, \mathcal{X})$ be a partition-regular pattern of complexity 1 such that $\mathcal{L}$ consists of $m$ linear forms. There exists a $\delta = \delta(\epsilon, m)$ with the following property. \cmnt{For all sufficiently large $n$, i}f $\phi: \vecsp{p}{n} \rightarrow [r]$ is an $r$-colouring of $\vecsp{p}{n}$ such that $\Lambda_{\mathcal{H}}(\phi) \leq \delta$, then there is a recolouring $\phi':\vecsp{p}{n} \rightarrow [r]$ that differs from $\phi$ in at most $\epsilon \abs{\vecsp{p}{n}}$ places such that $\Lambda_{\mathcal{H}}(\phi') = 0$.
\end{theorem}

\begin{proof}
    Let $A = \mathcal{K}(\mathcal{L})$ be an $\cmnt{\ell \times} m$ partition-regular matrix so that $\vecline{y} \in (\vecsp{p}{n})^m$ is an instance of $\mathcal{L}$ if and only if $A \underline{y} = 0$. Set $\delta' = (\epsilon/8r)^m/2m$ and $\epsilon' = \min(\epsilon, \epsilon_{count}(\delta'))$, where $\epsilon_{count}(\delta')$ is as in Definition \ref{def:complexity}. \cmnt{For $n$ sufficiently large, we may a}pply Lemma \ref{lemma:recolouring1} with $\epsilon'$ to obtain $\phi'$, three subspaces $H_3 \leqslant H_2 \leqslant H_1$ with codimensions $\codim(H_i) = D_i \leq C_{rcl}(\epsilon', r, m)$, a complement $U$ such that $U \oplus H_1 = \vecsp{p}{n}$, and $z_1, \ldots, z_m \in H_2$ satisfying $A\underline{z} = 0$.

    Suppose there exists an instance of $\mathcal{H}$ under $\phi'$, i.e.~$y_1, \ldots, y_m \in \vecsp{p}{n}$ such that $A(y_1, \ldots, y_m)^T = 0$ and $\chi \in \mathcal{X}$ such that $\phi'(y_i) = \chi(i)$ for all $i \in [m]$. Let $u_1, \ldots, u_m \in U$ be such that $y_i \in H_1 + u_i$.

    \begin{nclaim}
        The cosets $H_3+u_1+z_1, \ldots, H_3+u_m+z_m$ are consistent with $\mathcal{L}$.
    \end{nclaim}
    \begin{proof}[Proof of Claim]
    \renewcommand{\qedsymbol}{}
        For each $i$, write $y_i = h_i + u_i$ for some $h_1, \ldots, h_m \in H_1$. We know that $A(y_1, \ldots, y_m)^T = 0$ so
        $A(u_1, \ldots, u_m)^T = -A(h_1, \ldots, h_m)^T.$
        Since $H_1$ and $U$ are subspaces, $A\underline{h} \in H_1$ and $A\underline{u} \in U$, meaning that $A\underline{u} \in H_1 \cap U = \set{0}$. Therefore $(u_1, \ldots, u_m)$ is itself a solution to $A$. Since $(z_1, \ldots, z_m)$ is also a solution, so is $(u_1+z_1, \ldots, u_m+z_m)$, which completes the proof of the claim.
    \end{proof}

    \noindent As a shorthand, write $B_i$ for the coset $H_3 + u_i + z_i$. By Lemma \ref{lemma:recolouring1}\ref{prop:recolouring-regularity}, $\phi_{\chi(i)}$ is $\epsilon'$-uniform on $B_i$ for all $i \in [m]$, and $B_1, \ldots, B_m$ are consistent with $\mathcal{L}$ by the claim. Then the counting lemma (Lemma \ref{lemma:counting-lemma}) applies to give
    $$\Big|\Lambda_{\mathcal{H}}(\phi_{\chi_1} \id{B_1}, \ldots, \phi_{\chi_m} \id{B_m}) - p^{-D_3(m-\rank(\mathcal{L}))}\prod_{i=1}^m \alpha_i\Big|\leq p^{-D_3(m-\rank(\mathcal{L}))} m \delta',$$
    where $\alpha_i = \expct_{x \in B_i} \phi_{\chi(i)}(x) \geq \epsilon/8r$ by Lemma \ref{lemma:recolouring1}\ref{prop:recolouring-density}.

    Finally, note that $\Lambda_{\mathcal{H}}(\phi) \geq \Lambda_{\mathcal{H}}(\phi_{\chi_1} \id{B_1}, \ldots, \phi_{\chi_m} \id{B_m})$, so that, in particular, the density of $\mathcal{H}$ under $\phi$ is at least
    \vspace{-0.2cm}
    \begin{align*}
        \Lambda_{\mathcal{H}}(\phi_{\chi_1} \id{B_1}, \ldots, \phi_{\chi_m} \id{B_m}) &\geq p^{-D_3(m-\rank(\mathcal{L}))} \Big(\prod_{i=1}^m \alpha_i - m\delta'\Big)\\
        &\geq p^{-m C_{rcl}(\epsilon', r, m)}\Big((\epsilon/8r)^m - \frac12(\epsilon/8r)^m\Big)\\
        &\geq \frac12 \Big(8rp^{C_{rcl}(\epsilon', r, m)}\epsilon^{-1}\Big)^{-m},
    \end{align*}
    which contradicts the assumption of the theorem when $\delta < \left(8rp^{C_{rcl}(\epsilon', r, m)}\epsilon^{-1}\right)^{-m}/2$.
\end{proof}

\section{Discussion}
The largest contribution to the bound on $\delta$ in Theorem \ref{theorem:main-linear} is $C_{rcl}(\epsilon', r, m)$ arising from the use of Lemma \ref{lemma:recolouring1}. In turn, the proof of Lemma \ref{lemma:recolouring1} shows that $C_{rcl}(\epsilon', r, m)$ has the form $C_{rcl} = C_{reg}(\epsilon'p^{-C_{zreg}}/{8r}, \epsilon'/{4r}, \cmnt{r}, \ceil{\log_p({\epsilon}/{4r})}) + C_{zreg}$, where $C_{zreg} = C_{zreg}(\epsilon', \cmnt{r}, m)$ comes from Lemma \ref{lemma:third-partition} and $C_{reg}$ from Proposition \ref{prop:two-partitions}.

Lemma \ref{lemma:third-partition} follows from a single application of the arithmetic regularity lemma (Theorem \ref{theorem:arl}) and thus has the same growth in $\epsilon'$ as $C_{arl}$ in Theorem \ref{theorem:arl}. This is known to be of tower-type in general, with $\twr(n)$ defined recursively by $\twr(1) = 2$, $\twr(n) = 2^{\twr(n-1)}$. The lower bound of tower-type on $C_{arl}(\epsilon, 1)$ was first proved by Green in \cite{tower-type} as $\twr(\ceil{\log\epsilon^{-1}})$ and later strengthened by Hosseini et al.~in \cite{arl-lowerbound} to $\twr(\ceil{\epsilon^{-1}})$. As such, $C_{zreg}$ is $\twr(\ceil{1/\epsilon'})$ at best.

Proposition \ref{prop:two-partitions} follows from an iterated version of the arithmetic regularity lemma, known as the strong arithmetic regularity lemma (\cite[Theorem 5.4]{induced-1}). As a consequence of the iteration, one might expect wowzer-type growth in its bounds, where wowzer is an iterated tower defined by $\wwz(1) = 2$, $\wwz(n) = \twr(\wwz(n-1))$. Indeed, \cmnt{recent} work of the author \cite{sarl-lower-bound} confirms that wowzer-type growth in the strong arithmetic regularity lemma is unavoidable. This leads to wowzer-type bounds on $C_{reg}$ when the strong arithmetic regularity lemma is used.

However, Fox, Tidor and Zhao found an alternative proof of Proposition $\ref{prop:two-partitions}$ that bypasses the strong regularity lemma, resulting in better, tower-type bounds \cite[Section 5.2]{induced-1}. As a result, $C_{reg}(\epsilon' p^{-C_{zreg}}/{8r}, \epsilon/{4r}, \cmnt{r}, \ceil{\log_p({\epsilon}/{4r})})$ grows like a double tower in $1/\epsilon'$.

Finally, $\epsilon'$ was chosen to be at most $\epsilon_{count}(O(\epsilon^m))$ where $\epsilon_{count}$ is as in the definition of true complexity (Definition \ref{def:complexity}). In recent work, Manners \cite[Theorem 1.1.5]{fmanners} showed that $\epsilon_{count}(\delta)$ can be taken as polynomial in $\delta$.

\begin{theorem}[Polynomial bounds for true complexity \cite{fmanners}]
\label{thm:complexity-poly-bounds}
Let $\mathcal{L}$ be a linear system of finite complexity. Then for all $\delta > 0$, $\epsilon_{\ref{def:complexity}}(\delta)$ and $\epsilon_{count}(\delta)$ in Definition \ref{def:complexity} may be taken as $\delta^{O_{\mathcal{L}}(1)}$.
\end{theorem}
\noindent As a result, $\epsilon' = \epsilon^{O(1)}$, meaning that $C_{rcl}$, and therefore $\delta$, has growth of the order of $\twr(\twr(\epsilon^{-O(1)}))$.

It is not hard to modify the proof of Theorem \ref{theorem:main-linear} in order to remove any finite number of partition-regular patterns $\mathcal{H}_1 = (\mathcal{L}_1, \mathcal{X}_1), \ldots, \mathcal{H}_t = (\mathcal{L}_t, \mathcal{X}_t)$ at once (indeed, the induced removal lemmas of \cite{trans-invariant} and \cite{induced-1} are both able to remove multiple arithmetic patterns). This can be achieved by defining a single pattern encompassing all of $\mathcal{H}_1, \ldots, \mathcal{H}_t$. Let $A_1, \ldots, A_t$ denote the matrices $A_i = \mathcal{M}(\mathcal{L}_i)$, each of dimensions $\cmnt{\ell}_i \times m_i$. By taking $m = \max_i(m_i)$ and $\cmnt{\ell} = \max_i(\cmnt{\ell}_i)$, as well as adding zero rows or columns where appropriate, we may think of these as having the same dimension $\cmnt{\ell \times} m$ -- note that adding zero rows or columns does not affect partition regularity.

Let $A^*$ be a $t \ell \times tm$ block-diagonal matrix composed of $A_1, \ldots, A_t$, meaning that $A^*$ has the form
$$A^* = \begin{pmatrix}
    A_1 & \underline{0} & \ldots & \underline{0}\\
    \underline{0} & A_2 & \ldots & \underline{0}\\
     & & \ldots &\\
     \underline{0} & \underline{0} & \ldots & A_t
\end{pmatrix}.$$
Observe that $A^*$ inherits Rado's column conditions (Definition \ref{definition:column-condition}) from $A_1, \ldots, A_t$ via taking unions of the corresponding sets of columns. Thus, $A^*$ is itself partition-regular. Moreover, a solution to $A^* \underline{z} = 0$ is a vector $\underline{z} = (\underline{z}_1, \ldots, \underline{z}_t) \in (\vecsp{p}{n})^{mt}$ such that $A_i \underline{z}_i = 0$, i.e.~$\underline{z}_i$ is an instance of $\mathcal{L}_i$.

Additionally, given colourings $\chi_1, \ldots, \chi_t: [m] \rightarrow [r]$, let $(\chi_1, \ldots, \chi_t): [tm] \rightarrow [r]$ denote the joint colouring obtained by setting $(\chi_1, \ldots, \chi_t)(a) = \chi_{t'}(j)$
where $1 \leq t' \leq t$ and $1 \leq j \leq m$ are the unique integers satisfying $a = (t'-1)m + j$.

Let $\Psi = \set{\psi:[m] \rightarrow [r]}$ be the set of all possible $r$-colourings of $m$ and define $\mathcal{X}^*$ by 
$$\mathcal{X}^* = \set{(\chi_1, \ldots, \chi_t) \in \Psi^t \st \exists i \text{ s.t.~} \chi_i \in \mathcal{X}_i}.$$
Consider a new pattern $\mathcal{H}^* = (\mathcal{L}^*, \mathcal{X}^*)$. As noted, $\mathcal{H}^*$ is still partition-regular. Moreover, by definition, $\underline{z} = (\underline{z}_1, \ldots, \underline{z}_t) \in (\vecsp{p}{n})^{mt}$ is an instance of $\mathcal{H}^*$ if and only if each $\underline{z}_i$ is an instance of $\mathcal{L}_i$ and \textit{at least one of $\underline{z}_i$ is an instance of $\mathcal{H}_i$}. Therefore, if $\Lambda(\mathcal{H}_i)(\phi) \leq \delta$ for each $i$, then $\Lambda(\mathcal{H}^*)(\phi) \leq t\delta$. In particular, Theorem \ref{theorem:main-linear} applies to $\mathcal{H}^*$ to give a $\mathcal{H}^*$-recolouring. In fact, such a recolouring must also be $\mathcal{H}_i$-free for every $i \in [t]$, since for any instance $\underline{z}_i \in (\vecsp{p}{n})^m$ of $\mathcal{H}_i$, $(\underline{0}, \ldots, \underline{z}_i, \ldots, \underline{0})$ with $\underline{z}_i$ in the $i$th position is an instance of $\mathcal{H}^*$.

In fact, given this, it is possible to remove an infinite number of partition-regular patterns by reducing to the case of finitely many patterns, as was done, for instance, in \cite[Section 6]{induced-1} or \cite[Section 5.2]{trans-invariant}. 

Finally, it is worth noting that all arguments in Section \ref{section:subcosets} can be readily translated to apply to patterns of higher complexities, with the exception of a suitable version of Lemma \ref{lemma:third-partition}. In the higher-order setting, the types of partitions one obtains from the arithmetic regularity lemma are the level sets of polynomials (in the complexity 1 case, these polynomials are linear and thus result in cosets of a subspace). In order for such level sets to form a (nearly) equipartition, which is hard to forgo while counting, the polynomials need to be chosen carefully of high enough `rank' relative to each other (see \cite{montreal} for an excellent introduction to the quadratic case of this). To make the proof work, one would need a version of Lemma \ref{lemma:third-partition} that can be applied to the zero level set of the intermediate partition (the equivalent of $\mathcal{P}(H_2)$) in such a way that the polynomials of the new partition (the equivalent of $\mathcal{P}(H_3)$) interact well with the previously-obtained polynomials in terms of rank. However, this is hard to arrange without introducing a circular dependence between the various parameters in Lemma \ref{lemma:recolouring1}. It is the belief of the author that Theorem \ref{main} holds for partition-regular patterns of all finite complexities (and perhaps infinite complexity, as is the case in the induced removal lemma of Tidor and Zhao \cite{full-induced}), but the proof of this might require a novel way of handling the rank of partition refinements, or a different approach entirely.

\bibliographystyle{elsarticle-num} 
\bibliography{references}

\appendix
\renewcommand{\thesection}{\Alph{section}} % corrected redefinition of '\thesection'
\makeatletter
\renewcommand\@seccntformat[1]{\appendixname\ \csname the#1\endcsname.\hspace{0.5em}}
\makeatother

\section{Proof of the counting lemma}
%\section{Proof of the counting lemma}
\label{appendix:counting-lemma}

This appendix details the proof of the counting lemma stated in the introduction. Recall that $\epsilon_{count}$ was defined alongside true complexity in Definition \ref{def:complexity}.

\begin{replemma}{lemma:counting-lemma}
Fix $\delta > 0$, $\epsilon \leq \epsilon_{count}(\delta)$, and a linear system $\mathcal{L} = (L_1, \ldots, L_m)$ of complexity 1. Let $H \leqslant \vecsp{p}{n}$ be a subspace of codimension $d$. For any functions $f_1, \ldots, f_m:\vecsp{p}{n} \rightarrow [-1,1]$ and any $c_1, \ldots, c_m \in \vecsp{p}{n}$, if cosets $H+c_1, \ldots, H+c_m$ are consistent with $\mathcal{L}$ and $\epsilon$-regular for $f_1, \ldots, f_m$, then
$$\Big|\Lambda_{\mathcal{L}}(f_1 \id{H+c_1}, \ldots, f_m \id{H+c_m}) - p^{-d(m-\rank(\mathcal{L}))}\prod_{i=1}^m \alpha_i\Big| \leq p^{-d(m-\rank(\mathcal{L}))} m \delta,$$
where $\alpha_i$ denotes the average of $f_i$ on $H+c_i$.
\end{replemma}

\noindent We will first prove a short technical result that allows the density of the instances of $\mathcal{L}$ to be rewritten in terms of the matrix $\mathcal{K}(\mathcal{L})$ (see Section \ref{section:patterns}).
\begin{prop}
    \label{prop:appendix-counting-lemma-matrix}
    Let $\mathcal{L} = (L_1, \ldots, L_m)$ be a linear system in $l$ variables. If $A = \mathcal{K}(\mathcal{L})$, then for any functions $h_1, \ldots, h_m:\vecsp{p}{n} \rightarrow [-1,1]$,
    $\Lambda_{\mathcal{L}}(h_1, \ldots, h_m) = \expct_{\underline{z} \in \ker(A)} \prod_{i=1}^m h_i(z_i).$
\end{prop}
\begin{proof}
    Recall from Section \ref{section:patterns} that $M = \mathcal{M}(\mathcal{L})$ is the $m \times l$ matrix containing the coefficients of $L_i$ in the $i$th row, so that for any $\underline{x} \in (\vecsp{p}{n})^l$, $M\underline{x} = (L_1(\underline{x}), \ldots, L_m(\underline{x}))$. This means that $\Img{M}$ is precisely the same as the set of instances of $\mathcal{L}$. In particular,
    \begin{equation*}
        \sum_{\underline{x} \in (\vecsp{p}{n})^l} \prod_{i=1}^m h_i(L_i(\underline{x})) = \sum_{\underline{z} \in \Img{M}} \abs{\ker(M)} \prod_{i=1}^m h_i(z_i).
    \end{equation*}
    Moreover, $\ker(A) = \Img{M}$ by Proposition \ref{prop:solution-matrix} and $\abs{\vecsp{p}{n}}^l = \abs{\Img{M}}\abs{\ker(M)}$ by standard linear algebra, so
        $\Lambda_{\mathcal{L}}(h_1, \ldots, h_m) = \expct_{\underline{x} \in (\vecsp{p}{n})^l} \prod_{i=1}^m h_i(L_i(\underline{x})) = \expct_{\underline{z} \in \ker(A)} \prod_{i=1}^m h_i(z_i).$
\end{proof}

\begin{proof}[Proof of Lemma \ref{lemma:counting-lemma}]
    Let $A = \mathcal{K}(\mathcal{L})$ and write $\ker^H(A)$ for the restriction of $\ker(A)$ to $H^m$. Since $H+c_1,\ldots, H+c_m$ are consistent with $\mathcal{L}$, we may assume without loss of generality that $(c_1, \ldots, c_m)$ is itself an instance of $\mathcal{L}$, i.e.~$A\underline{c} = 0$. Then for any functions $h_1, \ldots, h_m:\vecsp{p}{n} \rightarrow [-1,1]$,
    \vspace{-0.1cm}
    $$\expct_{\underline{z} \in \ker(A)} \prod_{i=1}^m h_i(z_i)\id{H+c_i}(z_i) = \expct_{\underline{z} \in \ker(A)} \prod_{i=1}^m h_i(z_i)\id{H}(z_i-c_i) = \expct_{\underline{z} \in \ker^H(A)} \prod_{i=1}^m h_i(z_i+c_i),$$
    In particular, $\Lambda_{\mathcal{L}}(\id{H+c_1}, \ldots, \id{H+c_m}) = \abs{\ker^H(A)}/\abs{\vecsp{p}{n}} = p^{-d(m-\rank(\mathcal{L}))}$. Writing $\Lambda^H_{\mathcal{L}}$ for the density of $\mathcal{L}$-instances in $H$, we additionally have
    \begin{equation}
        \label{eq:appendix-counting-lemma-G-to-H}
        \Lambda_{\mathcal{L}}(f_1\id{H+c_1}, \ldots, f_m\id{H+c_m}) = p^{-d(m-\rank(\mathcal{L}))}\Lambda^H_{\mathcal{L}}(g_1, \ldots, g_m),
    \end{equation}
    where $g_i:H \rightarrow [-1,1]$ is the function given by $g_i(x) = f_i(x+c_i)$.

    On the other hand, using $\prod_{i=1}^m \alpha_i = \Lambda^H_{\mathcal{L}}(\alpha_1, \ldots, \alpha_m)$ in the telescoping identity \eqref{eq:telescoping} results in
    \begin{equation}
        \label{eq:appendix-counting-lemma-bound-eq}
        \Lambda^H_{\mathcal{L}}(g_1, \ldots, g_m) - \prod_{i=1}^m \alpha_i = \sum_{i=1}^m \Lambda^H_{\mathcal{L}}(h_1^\bracketed{i}, \ldots, h_m^\bracketed{i}),
    \end{equation}
    where for each $i$, at least one $1 \leq j \leq m$ satisfies $h^\bracketed{i}_{j} = g_j - \alpha_j$. Crucially, all $g_1, \ldots, g_m$ are $\epsilon$-uniform on $H$ by assumption, so that by Definition \ref{def:complexity}, $\babs{\Lambda^H_{\mathcal{L}}(h_1^\bracketed{i}, \ldots, h_m^\bracketed{i})} \leq \delta$. As a result, \eqref{eq:appendix-counting-lemma-bound-eq} can be bounded by $\babs{\Lambda^H_{\mathcal{L}}(g_1, \ldots, g_m) - \prod_{i=1}^m \alpha_i} \leq m\delta$, which together with \eqref{eq:appendix-counting-lemma-G-to-H} completes the proof.
\end{proof}

\end{document}